\documentclass[a4paper,12pt]{article}
\usepackage{amsmath}
\usepackage{amsfonts}
\usepackage{amssymb}
\usepackage{amsthm}

\newcommand{\al}{\alpha}
\newcommand{\be}{\beta}

\evensidemargin 0in \topskip .4in \theoremstyle{plain}
\newtheorem{thm}{Theorem}[section]

\newtheorem{lem}[thm]{Lemma}
\newtheorem{prop}[thm]{Proposition}
\theoremstyle{definition}
\newtheorem{defn}{Definition}[section]
\theoremstyle{remark}
\newtheorem{rem}{Remark}[section]
\numberwithin{equation}{section}

\begin{document}

\author{\textbf{A. S. Hegazi}$^{1}$\textbf{, W.Morsi}\textbf{, and M.Mansour}$^{2}$ \\
Mathematics department, Faculty of Science, \\
Mansoura university, Mansoura, 35516, EGYPT.\\
$^{1}$hegazi@mans.edu.eg,\\
$^{2}$mansour@mans.edu.eg}
\title{\textbf{Differential calculus on Hopf Group Coalgebra  }}

\date{.}
\maketitle
\begin{abstract}
In this paper we construct the Differential calculus on the Hopf
Group Coalgebra introduced by Turaev [10]. We proved that the
concepts introduced by S.L.Woronowicz in constructing Differential
calculus on Hopf Compact Matrix Pseudogroups (Quantum Groups)[7]
can be adapted to serve again in our construction.
\end{abstract}
\section*{Introduction}
Quantum groups, from a mathematical point of view, may be introduced by making emphasis on their $q-$deformed enveloping algebra aspects [1,2], which leads to the quantized enveloping algebras, or by making emphasis in the $R-$matrix formalism that describes the deformed group algebra. Also, they are mathematically well defined in the framework of Hopf algebra [3]. Quantum groups provide an interesting example of non-commutative geometry[4]. Non-commutative differential calculus on quantum groups is a fundamental tool needed for many applications [5,6].\\
\indent
S.L.Woronowicz [7] gave the general framework for bicovariant differential calculus on quantum groups following general ideas of A.Connes. Also, He showed that all important notions and formulae of classical Lie group theory admit a generalization to the quantum group case and he has restricted himself to compact matrix pseudogroups as introduced in [8].In contrast to the classical differential geometry on Lie groups, there is no functorial method to obtain a unique bicovariant differential calculus on a given quantum group [9]. \\
\indent
Recently, Quasitriangular Hopf $\pi-$coalgebras are introduced by Turaev [10]. He has showed that they give rise to crossed $\pi-$categories. Virelizier [11] studied the algebraic properties of the Hopf $\pi-$coalgebras, also he has showed that the existence of integrals and trace for such coalgebras and has generalized the main properties of the quasitriangular Hopf algebras to the setting of Hopf $\pi-$coalgebra.\\
\indent
In this paper we will use the concepts introduced by S.L.Woronowicz [7] to construct the Differential calculus on the Hopf group coalgebra(introduced by Turaev [10]). We briefly describe the content of the paper.In section one we give the definition of Hopf group coalgebras [11]. In section two, we give the main definitions and theorems concerning first order differential calculus. Section three contains the construction of the $\pi-$graded Bicovariant bimodules. Finally, in section four we construct the first order differential calculus on the Hopf group coalgebra.\\
Now let us give some  basic definitions about Hopf $\pi-$coalgebra \\
\section{Hopf Group Coalgebra}
\begin{defn}
 A $\pi-$coalgebra is a family $C=\{C_{\alpha}\}_{\alpha\in \pi}$ of $\ k-$linear
spaces endowed with a family
$\Delta=\{\Delta_{\alpha,\beta}:C_{\alpha\beta} \to C_{\alpha}
\otimes C_{\beta}\}_{\alpha,\beta\in \pi}$ of $\ k-$linear maps
(the comultiplication) and a $\ k-$linear map $\varepsilon
:C_{1}\to\ k$
 such that
\begin{itemize}
    \item $\Delta$ is coassociative in the sense that for any $\alpha,\beta,\gamma\in\pi$,
$$
    (\Delta_{\alpha,\beta}\otimes id)\Delta_{\alpha\beta,\gamma}=(id\otimes\Delta_{\beta,\gamma})\Delta_{\alpha,\beta\gamma}\;,
$$
    \item for all $\alpha \in \pi$,
    $$
    (id\otimes\varepsilon)\Delta_{\alpha, 1}=(\varepsilon \otimes id)\Delta_{1,\alpha} .
    $$
\end{itemize}
\end{defn}
\noindent \textbf{Sweedler's notation} In the case of Hopf group
 coalgebra Sweedler's notations have been extended by Turaev and Virelizier
 in the following way: for any $\alpha,\beta\in \pi$ and $c\in C_{\alpha\beta}$,
 they defined
$$
\Delta_{\alpha,\beta}(c)=\sum_{(c)}c_{(1,\alpha)} \otimes
c_{(2,\beta)}\in C_{\alpha}\otimes C_{\beta}.
$$
or shortly, if we have  the summation implicit
$$
\Delta_{\alpha,\beta}(c)=c_{(1,\alpha)} \otimes c_{(2,\beta)}.
$$
The coassociativity axiom gives that , for any
$\alpha,\beta,\gamma\in \pi $ and $c\in C_{\alpha\beta,\gamma}$
$$
c_{(1,\alpha\beta)(1,\alpha)}\otimes
c_{(1,\alpha\beta)(2,\beta)}\otimes
c_{(2,\gamma)}=c_{(1,\alpha)}\otimes
c_{(2,\beta\gamma)(1,\beta)}\otimes c_{(2,\beta\gamma)(2,\gamma)}
.$$

\medskip

\indent Let $C=(\{C_{\alpha}\}_{\alpha
 \in \pi},\Delta,\varepsilon)$ be a
$\pi-$coalgebra and $A$ be an algebra with multiplication $m$ and
unit element $1_{A}$. The family $\Delta$ and the map $m$ induce a
map
$$
*:conv(C,A)\otimes conv(C,A)\to conv(C,A)
$$
defined by the composition
\begin{eqnarray*}
Hom(C_{\alpha},A)\otimes
Hom(C_{\beta},A)\stackrel{\rho}{\longrightarrow}&Hom(C_{\alpha}\otimes
C_{\beta},A\otimes
A)&\stackrel{Hom(\Delta_{\alpha,\beta},m)}{\longrightarrow}\\&&Hom(C_{\alpha
\beta},A)
\end{eqnarray*}
where $\rho$ is the natural injection of $Hom(C_{\alpha},A)\otimes
Hom(C_{\beta},A)$ into $Hom(C_{\alpha}\otimes C_{\beta},A\otimes
\nolinebreak A)$\\
The map $*$ is called convolution product of $f,g$\\
Also, the maps
$$ \varepsilon :C_{1}\longrightarrow  \ k \;\;\;\;\textrm{and} \;\;\;\;
\eta:\ k \longrightarrow A$$ induce a map
$$ \eta_{Conv(C,A)}:\ k\longrightarrow Conv(C,A)$$
defined by
$$ (\eta_{Conv(C,A)}(\lambda))(c)=\varepsilon (c)\eta(\lambda)$$
for all $c\in C_{1}$.

\begin{lem}
The $\ k$-space $$ Conv(C,A) =\bigoplus_{\al\in\pi}
Hom(C_{\al},A)$$ endowed with the convolution product $*$ and the
unit element $\varepsilon 1_{A}$ is a $\pi-$graded algebra called
the convolution algebra.
\end{lem}

\begin{rem}
If we put $A=\ k$ in the above lemma the $\pi-$graded algebra
$Conv(C,\ k)=\bigoplus_{\alpha\in \pi}C^{*}_{\alpha}$ is called
dual to $C$ and denoted by $C^{*}$.
\end{rem}

\begin{defn}
 A Hopf $\pi-$ coalgebra is a $\pi-$coalgebra  $H=(\{H_{\alpha} \}_{\alpha
 \in \pi}, \Delta,\varepsilon )$ endowed with a family
$$
S=\{S_{\alpha}:H_{\alpha} \to H_{\alpha^{-1}}\}_{\alpha\in\pi}$$
of $\ k-$linear maps called the antipode such that

\begin{enumerate}
  \item[(1)] Each $H_{\alpha}$ is an algebra with multiplication $m_{\alpha}$
  and unit element $1_{\alpha}\in H_{\alpha}$,
  \item[(2)] The linear maps
  \begin{eqnarray*}
  \Delta_{\alpha,\beta}&:&H_{\alpha\beta}\to H_{\alpha}\otimes H_{\beta},\\
  \varepsilon &:& H_{1}\to \ k.
  \end{eqnarray*}
  are algebra maps for all $\alpha,\beta\in A$,
  \item[(3)] For any $\alpha\in \pi$
  $$m_{\alpha}(S_{\alpha^{-1}}\otimes id)\Delta_{\alpha^{-1},\alpha}=
  m_{\alpha}(id\otimes S_{\alpha^{-1}})\Delta_{\alpha,\alpha^{-1}}.$$
\end{enumerate}
\end{defn}

\medskip

\begin{rem}
If $H=(\{H_{\alpha}\}_{\alpha\in \pi},\Delta,\varepsilon ,S)$ is a
Hopf $\pi-$coalgebra then axiom (3) says that $S_{\alpha}$ is the
inverse of $I_{H_{\alpha^{-1}}}$ in the convolution algebra
$Conv(H,H_{\alpha^{-1}})$.
\end{rem}

\medskip

\begin{rem}
$(H_{1},\Delta_{1,1},\varepsilon ,S_{1} )$ is a classical Hopf
algebra
\end{rem}
\begin{lem}
Let $H=(\{H_{\alpha}\}_{\alpha\in \pi},\Delta,\varepsilon ,S)$ be
a Hopf $\pi-$coalgebra. then
\begin{enumerate}
\item
$\Delta_{\beta^{-1},\alpha^{-1}}S_{\alpha\be}=\sigma_{H^{\alpha^{-1}}
,H^{\beta^{-1}}}(S_{\alpha}\otimes S_{\beta}
)\Delta_{\alpha,\beta}$ for any $\alpha,\beta\in \pi$,
  \item  $\varepsilon(S_{1})=\varepsilon$,
  \item $S_{\alpha} (ab)= S_{\alpha} (b)S_{\alpha} (a)$ for any $\alpha\in\pi$ and $a,b\in A$,
  \item $S_{1_{\alpha}} =1_{\alpha^{-1}}$ for any $\alpha\in\pi$.
\end{enumerate}
\end{lem}
\begin{defn}
Let $C=(\{C_{\alpha}\}_{\alpha\in\pi},\Delta,\varepsilon)$ be a
$\pi-$coalgebra. A right $\pi-$comodule over $C$ is a family
$M=\{M_{\alpha}\}_{\alpha\in\pi}$ of $\ k-$linear spaces endowed
with a family $\rho=\{\rho_{\alpha,\beta}:
M_{\alpha\beta\longrightarrow M_{\alpha}\otimes C_{\beta}}\}$ of
$\ k-$linear maps (the structure maps) such that
\begin{itemize}
    \item For any $\alpha,\beta,\gamma\in\pi$
    $$(\rho_{\alpha,\beta}\otimes id)\rho_{\alpha\beta,\gamma}=(id\otimes \Delta_{\beta,\gamma})\rho_{\alpha,\beta\gamma}\eqno{*}$$
    \item For any $\alpha\in\pi$
    $$
(id\otimes \varepsilon)\rho_{\alpha,1}=id \eqno{**}
    $$
\end{itemize}
\end{defn}

\begin{defn}
A $\pi-$subcomodule of $M$ is a family
$N=\{N_{\alpha}\}_{\alpha\in\pi}$ where $N_{\alpha}$ is a $\
k-$linear subspace of $M_{\alpha} $ such that for all
$\alpha,\beta\in \pi $
$$
\rho_{\alpha,\beta}(N_{\alpha\beta})\subset N_{\alpha}\otimes
C_{\beta}
$$
\end{defn}

\begin{defn}
A $\pi-$comodule morphism between to right $\pi-$comodules $M$ and
$M'$ over a $\pi-$coalgebra $C$ (with structure maps $\rho$ and
$\rho'$, respectively) is a family
$f=\{f_{\alpha}:M_{\alpha}\longrightarrow M_{\alpha}'\}$ of $\
k-$linear maps such that for all $\alpha,\beta\in\pi$
$$
\rho_{\alpha,\beta}'(f_{\alpha\beta})=(f_{\alpha}\otimes
id)\rho_{\alpha,\beta}
$$
\end{defn}

\noindent\textbf{Sweedler's notation}\\
For any $\alpha,\beta \in \pi$ and $m\in M_{\alpha,\beta}$ we
write

$$
 \rho_{\alpha,\beta}(m)=m_{(0,\alpha)}\otimes m_{(1,\beta)}\in
 M_{\alpha}\otimes C_{\beta}
$$
also the axiom
    $$(\rho_{\alpha,\beta}\otimes id)\rho_{\alpha\beta,\gamma}=(id\otimes
    \Delta_{\beta,\gamma})\rho_{\alpha,\beta\gamma}$$
can be written as
$$
m_{(0,\alpha\beta)(0,\alpha)}\otimes
m_{(0,\alpha\beta)(1,\beta)}\otimes m_{(1,\gamma)}
=m_{(0,\alpha)}\otimes m_{(1,\beta\gamma)(1,beta)}\otimes
m_{(1,\beta\gamma)(2,\gamma)}$$

This elements of $M_{\alpha}\otimes C_{\beta}\otimes C_{\gamma}$
is written as $m_{(0,\alpha)}\otimes m_{(1,\beta)}\otimes
m_{(2,\gamma)}$

\section{Basic Definitions of differential calculus}
\begin{defn}  Let $A=\{ A_{\alpha }\} _{\alpha \in \pi }$ be a Hopf group
coalgebra ,$\Gamma =\{ \Gamma _{\alpha }\} _{\alpha \in \pi }$ be
a $\pi -$ graded bimodule over $A$ , and
$$ d=\{ d_{\alpha }:A_{\alpha
}\longrightarrow \Gamma _{\alpha }\} \eqno{2.1}$$

 be a family of linear maps. We say that $\left( \Gamma ,d\right) $ is a $\pi -
$graded first order differential  calculus over $A$ if for any
$\alpha \in \pi $
\begin{enumerate}
    \item For any $a,b\in A_{\alpha }$
          $$d_{\alpha }( ab) =d_{\alpha }( a)
b+ad_{\alpha }( b)  \eqno{2.2}$$

    \item Any element $\rho \in \Gamma _{\alpha }$ is of the form

$$\rho =\sum_{k=1}^{n} a_{k}d_{\alpha }b_{k }\;\; ,\;\;a_{k},b_{k}\in A_{\alpha }$$
\end{enumerate}
\end{defn}
%==========================================================================================
\medskip
\begin{defn}
 Two $\pi -$graded first order differential calculi are said to
be isomorphic if there exists a family of bimodule isomorphisms
$i=\left\{ i_{\alpha }:\Gamma _{\alpha }\longrightarrow \Gamma
_{\alpha }^{^{\prime }}\right\} $ such that

$$i_{\alpha }\left( d_{\alpha }a\right) =d_{\alpha }^{^{\prime }}a
,\mbox{ for all } a\in A_{\alpha },\alpha \in \pi .$$
 \end{defn}
%===========================================================================================

Let $A=\left\{ A_{\alpha }\right\} _{\alpha \in \pi }$ be a Hopf
group coalgebra ,  $m_{\alpha }:A_{\alpha }\otimes A_{\alpha
}\longrightarrow A_{\alpha }$ be the multiplication defined on
$A_{\alpha }$ for each $\alpha $. Define $A^{2}$ $=\left\{
A_{\alpha }^{2}\right\} _{\alpha \in \pi }$ such that

$$A_{\alpha }^{2}=\left\{ q\in A_{\alpha }\otimes A_{\alpha },m_{\alpha
}\left( q\right) =0\right\} \eqno{2.3}$$

By definition $A_{\alpha }^{2}$ is a linear subspace of $A_{\alpha
}\otimes A_{\alpha }$ for each $\alpha \in \pi $ .On $A^{2}$
define an $A-$bimodule structure as

For any $\alpha \in \pi ,c\in $ $A_{\alpha },\sum_{k}a_{k}\otimes
b_{k}\in A_{\alpha }^{2}$

$$c\left( \sum_{k}a_{k}\otimes b_{k}\right) =\sum_{k}ca_{k}\otimes
b_{k}\eqno{2.4}$$

$$\left( \sum_{k}a_{k}\otimes b_{k}\right) c=\sum_{k}a_{k}\otimes
b_{k}c\eqno{2.5}$$

Define $D=\left\{ D_{\alpha }\right\} _{\alpha \in \pi }$ by
$$D_{\alpha }\left( b\right) =1_{\alpha }\otimes b-b\otimes
1_{\alpha },$$
 for all $ b\in A_{\alpha },\alpha \in \pi $
\\

It is clear that $m_{\alpha }\left( D_{\alpha }\left( b\right) \right) =0,$%
i.e.$D_{\alpha }\left( b\right) \in A_{\alpha }^{2}$. Moreover
\begin{eqnarray*}
D_{\alpha }\left( ab\right) &=&1_{\alpha }\otimes ab-ab\otimes
1_{\alpha }\\
&=&1_{\alpha }\otimes ab-a\otimes b+a\otimes b-ab\otimes 1_{\alpha
}\\
&=&\left( 1_{\alpha }\otimes a-a\otimes 1_{\alpha }\right)
b+a\left( 1_{\alpha }\otimes b-b\otimes 1_{\alpha }\right) \\
&=&D_{\alpha }\left( a\right) b+aD_{\alpha }\left( b\right)
\end{eqnarray*}
%==========================================================================================
\begin{prop}
Let $N=\left\{ N_{\alpha }\right\} _{\alpha \in \pi }$ be a $\pi
-$graded sub-bimodule of $A^{2},\Gamma =A^{2}/N$ , $\pi =\left\{
\pi _{\alpha }:A_{\alpha }^{2}\longrightarrow \Gamma _{\alpha
}\right\} $ be the family of canonical epimorphisms , and
$d=\left\{ d_{\alpha }=\pi _{\alpha }\circ D_{\alpha }\right\}
_{\alpha \in \pi }.$ Then $\Gamma =\left( \left\{ \Gamma_{\alpha}
\right\} _{\alpha \in \pi },d\right) $ is a first order
differential calculus over $A.$Any other $\pi -$graded first order
differential calculus over $A$ can be obtained in this way.
\end{prop}

\begin{proof}
By definition of $\Gamma =\left\{ \Gamma_{\alpha} \right\}
_{\alpha \in \pi },\Gamma $ is a $\pi -$graded bimodule over
$A.$Moreover ,by definition of $d=\left\{ d_{\alpha }=\pi _{\alpha
}\circ D_{\alpha }\right\} _{\alpha \in \pi }$ we find that
$\Gamma =\left( \left\{ \Gamma_{\alpha} \right\} _{\alpha \in \pi
},d\right) $ is a $\pi -$graded first order differential calculus
over $A.$ It remains to show that any $\pi -$graded first order
differential calculus over $A$ can be obtained in this way.\\
 Let $\Gamma =\left( \left\{ \Gamma _{\alpha }\right\} _{\alpha \in \pi
},d\right) $ be any other $\pi -$graded first order differential
calculus over $A$ .We have for each $\alpha \in \pi ,$
$\sum_{k}a_{k}\otimes b_{k}\in $ $A_{\alpha }^{2},c\in $
$A_{\alpha }$

$$ \sum_{k}ca_{k}d_{\alpha }b_{k}=c\left( \sum_{k}a_{k}d_{\alpha
}b_{k}\right) $$ and
\begin{eqnarray*}
\sum_{k}a_{k}d_{\alpha }\left( b_{k}c\right) &=&\left(
\sum_{k}a_{k}d_{\alpha }b_{k}\right) c+\left(
\sum_{k}a_{k}b_{k}\right) d_{\alpha }c\\
&=&\left( \sum_{k}a_{k}d_{\alpha }b_{k}\right) c
\end{eqnarray*}
i.e. the family $\pi =\left\{ \pi _{\alpha }:A_{\alpha
}^{2}\longrightarrow \Gamma _{\alpha }\right\} $defined by the
formula

$$\pi _{\alpha }\left( \sum_{k}a_{k}\otimes b_{k}\right) =\sum_{k}a_{k}d_{\alpha }b_{k}\eqno{2.6}$$
is a bimodule morphism.We will show that $\pi _{\alpha }$ is
surjective for each $\alpha \in \pi .$

Let $\rho \in \Gamma _{\alpha }$ such that $$\rho
=\sum_{k}a_{k}d_{\alpha }b_{k},\;\;\;\;\;a_{k},b_{k}\in A_{\alpha
}$$ Define an element $q\in $ $A_{\alpha }\otimes $ $A_{\alpha }$
by

$$q=\sum_{k}a_{k}\otimes b_{k}-a_{k}b_{k}\otimes 1_{\alpha }$$
It is clear that $m_{\alpha }q=0$ , i.e. $q$ $\in $ $A_{\alpha
}^{2}$ .Moreover ,
\begin{eqnarray*}
\pi _{\alpha }\left( q\right) &=&\sum_{k}a_{k}d_{\alpha
}b_{k}-a_{k}b_{k}d_{\alpha }1_{\alpha }\\
&=&\sum_{k}a_{k}d_{\alpha }b_{k}\\
 &=&\rho
\end{eqnarray*}
therefore $\pi _{\alpha }$ is surjective for each $\alpha \in
\pi$.
\begin{eqnarray*}
\ker \pi &=&\left\{ \ker \pi _{\alpha }\right\} _{\alpha \in \pi }\\
&=&\left\{ \sum_{k}a_{k}\otimes b_{k}\in A_{\alpha }^{2}
 , \sum_{k}a_{k}d_{\alpha }b_{k}=0\right\} _{\alpha \in \pi }
\end{eqnarray*}
Taking
$$N=\{ N_{\alpha }=\ker \pi _{\alpha }=\{ \sum_{k}a_{k}\otimes
b_{k}\in A_{\alpha }^{2}\sum_{k}a_{k}d_{\alpha }b_{k}=0\}
\}_{\alpha \in \pi }\eqno{2.7}$$ then $\Gamma $ can be identified
by $A^{2}/N$ and for any $b\in A_{\alpha }$
\begin{eqnarray*}
\pi _{\alpha }D_{\alpha }\left( b\right) &=&\pi _{\alpha }\left(
1_{\alpha }\otimes b-b\otimes 1_{\alpha }\right)\\
&=&d_{\alpha }b-bd_{\alpha }1_{\alpha }\\
 &=&d_{\alpha }b.
\end{eqnarray*}
\end{proof}
%===========================================================================================
\begin{defn}

Let $\Gamma =\left( \left\{ \Gamma_{\alpha} \right\} _{\alpha \in
\pi },d\right) $ be a $\pi -$graded first order differential
calculus over $A$.We say that $ \Gamma =\left( \left\{
\Gamma_{\alpha} \right\} _{\alpha \in \pi },d\right) $ is left
covariant if for any $\alpha ,\beta \in \pi $

$$\sum_{k}a_{k}d_{\alpha \beta }b_{k} \Longrightarrow \sum_{k}\Delta _{\alpha ,\beta }\left( a_{k}\right) \left(
id\otimes d_{\beta }\right) \Delta _{\alpha ,\beta }\left(
b_{k}\right) =0\eqno{2.8}$$ for any $a_{k},b_{k}\in A_{\alpha
\beta },k=1,2,......,n.$
\end{defn}

\begin{prop}
Let $\Gamma =\left( \left\{ \Gamma_{\alpha} \right\} _{\alpha \in
\pi },d\right) $ be
a left covariant $\pi -$graded first order differential calculus over $A.$%
Then there exists a family of linear mappings

$$ \Delta ^{l}=\left\{ \Delta _{\alpha ,\beta }^{l}:\Gamma _{\alpha
\beta }\longrightarrow A_{\alpha }\otimes \Gamma _{\beta }\right\}
\eqno{2.9}$$

such that
\begin{enumerate}
    \item For any $a\in A_{\alpha \beta },\rho \in \Gamma _{\alpha
    \beta}$

$$\Delta _{\alpha ,\beta }^{l}\left( a\rho \right) =\Delta _{\alpha
,\beta }\left( a\right) \Delta _{\alpha ,\beta }^{l}\left( \rho
\right)\eqno{2.10}$$

$$\Delta _{\alpha ,\beta }^{l}\left( \rho a\right) =\Delta _{\alpha
,\beta }^{l}\left( \rho \right) \Delta _{\alpha ,\beta }\left(
a\right) \eqno{2.11}$$

    \item For any $\alpha,\beta,\gamma\in \pi$
    $$\left( \Delta _{\alpha ,\beta }\otimes id\right)
\Delta _{\alpha \beta ,\gamma }^{l}=\left( id\otimes \Delta
_{\beta ,\gamma }^{l}\right) \Delta _{\alpha ,\beta \gamma
}^{l}\eqno{2.12}$$

    \item For any $\rho \in \Gamma _{\alpha }$

$$\left( \varepsilon \otimes id\right) \Delta
_{1,\alpha }^{l}\left( \rho \right) =\rho \eqno{2.13}$$
    \item For any $\alpha,\beta \in \pi$

    $$\Delta _{\alpha ,\beta }^{l}d_{\alpha
\beta }=\left( id\otimes d_{\beta }\right) \Delta _{\alpha ,\beta
}\left( a\right)$$
\end{enumerate}
\end{prop}

\begin{proof}

Let $\Delta ^{l}=\left\{ \Delta^{l} _{\alpha ,\beta }\right\}
_{\alpha ,\beta \in \pi }$ where $\Delta _{\alpha ,\beta
}^{l}:\Gamma _{\alpha \beta }\longrightarrow A_{\alpha }\otimes
\Gamma _{\beta }$ is defined by
$$\Delta _{\alpha ,\beta }^{l}\left(
\sum_{k=1}^{n}a_{k}d_{\alpha \beta }b_{k}\right)
=\sum_{k=1}^{n}\Delta _{\alpha ,\beta }\left( a_{k}\right) \left(
id\otimes d_{\beta }\right) \Delta _{\alpha ,\beta }\left(
b_{k}\right) $$
\noindent where $a_{k},b_{k}\in A_{\alpha \beta
},\alpha ,\beta \in \pi$. Then by definition for each $\alpha
,\beta \in \pi $ $\Delta _{\alpha ,\beta }^{l}$ is a well defined
linear map.
\begin{enumerate}
\item Let $a\in A_{\alpha \beta }$ $,$ $\rho \in \Gamma _{\alpha
\beta },$ $\rho =\sum_{k=1}^{n}a_{k}d_{\alpha \beta
}b_{k},a_{k},b_{k}\in A_{\alpha \beta }$
\begin{eqnarray*}
\Delta _{\alpha ,\beta }^{l}\left( \rho a\right) &=&\Delta
_{\alpha ,\beta }^{l}\left( \left(
\sum_{k}a_{k}d_{\alpha \beta }b_{k}\right) a\right)\\
&=&\Delta _{\alpha ,\beta }^{l}\left( \sum_{k}a_{k}d_{\alpha \beta
}\left( b_{k}a\right) -\sum_{k}a_{k}b_{k}d_{\alpha \beta }a\right) \\
&=&\sum_{k}\left( \Delta _{\alpha ,\beta }\left( a_{k}\right)
\left( id \otimes d_{\beta }\right) \Delta _{\alpha ,\beta }\left(
b_{k}a\right) -\Delta _{\alpha ,\beta }\left( a_{k}b_{k}\right)
\left( id\otimes d_{\beta }\right)
\Delta _{\alpha ,\beta }\left( a\right) \right)\\
&=&\sum_{k}\Delta _{\alpha ,\beta }\left( a_{k}\right) \left(
\left( id\otimes d_{\beta }\right) \Delta _{\alpha ,\beta }\left(
b_{k}a\right) -\Delta _{\alpha ,\beta }\left( b_{k}\right) \left(
id\otimes d_{\beta }\right) \Delta
_{\alpha ,\beta }\left( a\right) \right)\\
&=&\sum_{k}\Delta _{\alpha ,\beta }\left( a_{k}\right) (\left(
id\otimes d_{\beta }\right) \left( b_{k\left( 1,\alpha \right)
}a_{\left( 1,\alpha \right) }\otimes b_{k\left( 2,\beta \right)
}a_{\left( 2,\beta \right) }\right) -b_{k\left( 1,\alpha \right)
}a_{\left( 1,\alpha \right) }\otimes b_{k\left( 2,\beta
\right) }d_{\beta }a_{\left( 2,\beta \right) })\\
&=&\sum_{k}\Delta _{\alpha ,\beta }\left( a_{k}\right) (b_{k\left(
1,\alpha \right) }a_{\left( 1,\alpha \right) }\otimes d_{\beta
}\left( b_{k\left( 2,\beta \right) }a_{\left( 2,\beta \right)
}\right) -b_{k\left( 1,\alpha \right) }a_{\left( 1,\alpha \right)
}\otimes b_{k\left( 2,\beta \right) }d_{\beta }a_{\left( 2,\beta
\right) })\\
&=&\sum_{k}\Delta _{\alpha ,\beta }\left( a_{k}\right)
[(b_{k\left( 1,\alpha \right) }a_{\left( 1,\alpha \right) }\otimes
d_{\beta }b_{k\left( 2,\beta \right) }a_{\left( 2,\beta \right)
}+b_{k\left( 1,\alpha \right) }a_{\left( 1,\alpha \right) }\otimes
b_{k\left( 2,\beta \right) }d_{\beta }a_{\left( 2,\beta \right)
})-b_{k\left( 1,\alpha \right) }a_{\left( 1,\alpha \right)
}\\&&\otimes b_{k\left( 2,\beta \right) }d_{\beta }a_{\left(
2,\beta
\right) }]\\
&=&\sum_{k}\Delta _{\alpha ,\beta }\left( a_{k}\right) \left(
b_{k\left( 1,\alpha \right) }a_{\left( 1,\alpha \right) }\otimes
d_{\beta }b_{k\left( 2,\beta \right) }a_{\left( 2,\beta \right)
}\right) \\
&=&\left( \sum_{k}\Delta _{\alpha ,\beta }\left( a_{k}\right)
\left( id\otimes d_{\beta }\right) \Delta _{\alpha ,\beta }\left(
b_{k}\right) \right) \Delta _{\alpha ,\beta }\left(
a\right)\\
&=&\Delta _{\alpha ,\beta }^{l}\left( \sum_{k}a_{k}d_{\alpha \beta
}b_{k}\right) \Delta _{\alpha ,\beta }\left( a\right)\\
&=&\Delta _{\alpha ,\beta }^{l}\left( \rho \right) \Delta _{\alpha
,\beta }\left( a\right)
\end{eqnarray*}
and
\begin{eqnarray*}
 \Delta_{\alpha ,\beta }^{l}\left( a\rho \right) &=&\Delta _{\alpha
,\beta }^{l}\left( a\sum_{k}a_{k}d_{\alpha \beta
}b_{k}\right)\\
&=&\Delta _{\alpha ,\beta }^{l}\left( \sum_{k}aa_{k}d_{\alpha
\beta }b_{k}\right)\\
&=&\sum_{k}\Delta _{\alpha ,\beta }\left( aa_{k}\right) \left(
id\otimes d_{\beta }\right) \Delta _{\alpha ,\beta
}\left( b_{k}\right)\\
&=&\Delta _{\alpha ,\beta }\left( a\right) \sum_{k}\Delta _{\alpha
,\beta }\left( a_{k}\right) \left( id\otimes d_{\beta
}\right) \Delta _{\alpha ,\beta }\left( b_{k}\right)\\
&=&\Delta _{\alpha ,\beta }\left( a\right) \Delta _{\alpha ,\beta
}^{l}\left( \sum_{k}a_{k}d_{\alpha \beta }b_{k}\right)\\
&=&\Delta _{\alpha ,\beta }\left( a\right) \Delta _{\alpha ,\beta
}^{l}\left( \rho \right)\\
\end{eqnarray*}
%=================================

\item  Let $ad_{\alpha \beta \gamma }b\in \Gamma _{\alpha \beta \gamma },$with $%
a,b\in A_{\alpha \beta \gamma },$ then we compute
\begin{eqnarray*}
\left( \Delta _{\alpha ,\beta }\otimes id\right) \Delta _{\alpha
\beta ,\gamma }^{l}\left( ad_{\alpha \beta \gamma }b\right)
&=&\left( \Delta _{\alpha ,\beta }\otimes id\right) \left( \Delta
_{\alpha \beta ,\gamma }\left( a\right) \left( id\otimes d_{\gamma
}\right) \Delta _{\alpha \beta ,\gamma }\left(
b\right) \right)\\
&=&\left( \Delta _{\alpha ,\beta }\otimes id\right) \Delta
_{\alpha \beta ,\gamma }\left( a\right) \left( \Delta _{\alpha
,\beta }\otimes id\right) \left( id\otimes d_{\gamma }\right)
\Delta _{\alpha \beta ,\gamma }\left(
b\right)\\
&=&\left( \Delta _{\alpha ,\beta }\otimes id\right) \Delta
_{\alpha \beta ,\gamma }\left( a\right) \left( id\otimes id\otimes
d_{\gamma }\right) \left( \Delta _{\alpha ,\beta }\otimes
id\right) \Delta _{\alpha \beta ,\gamma }\left( b\right)
\end{eqnarray*}
On the other hand
\begin{eqnarray*}
\left( id\otimes \Delta _{\beta ,\gamma }^{l}\right) \Delta
_{\alpha ,\beta \gamma }^{l}\left( ad_{\alpha \beta \gamma
}b\right)&=&\left( id\otimes \Delta _{\beta ,\gamma }^{l}\right)
\left( \Delta _{\alpha ,\beta \gamma }\left( a\right) \left(
id\otimes d_{\beta \gamma }\right) \Delta _{\alpha
,\beta \gamma }\left( b\right) \right)\\
&=&\left( id\otimes \Delta _{\beta ,\gamma }^{l}\right) \left(
a_{\left( 1,\alpha \right) }b_{\left( 1,\alpha \right) }\otimes
a_{\left( 2,\beta \gamma \right) }d_{\beta \gamma }b_{\left(
2,\beta \gamma \right) }\right)\\
&=&a_{\left( 1,\alpha \right) }b_{\left( 1,\alpha \right) }\otimes
a_{\left( 2,\beta \right) }b_{\left( 2,\beta \right) }\otimes
a_{\left( 3,\gamma \right) }d_{\gamma }b_{\left( 3,\gamma \right)
}\\
&=&\left( id\otimes \Delta _{\beta ,\gamma }\right) \Delta
_{\alpha ,\beta \gamma }\left( a\right) \left( id\otimes id\otimes
d_{\gamma }\right) \left( id\otimes \Delta _{\beta ,\gamma
}\right) \Delta _{\alpha ,\beta \gamma }\left(
b\right)\\
\end{eqnarray*}

\item For $\alpha \in \pi $ let $ad_{\alpha }b\in \Gamma _{\alpha
},$ $a,b\in A_{\alpha }$
\begin{eqnarray*}
\left( \varepsilon \otimes id\right) \Delta _{1,\alpha }^{l}\left(
ad_{\alpha }b\right) &=&\left( \varepsilon \otimes id\right)
\left( \Delta _{1,\alpha }\left( a\right) \left( id\otimes
d_{\alpha }\right) \Delta _{1,\alpha }\left( b\right) \right)\\
&=&\varepsilon \left( a_{\left( 1,1\right) }b_{\left( 1,1\right)
}\right) a_{\left( 2,\alpha \right) }d_{\alpha }b_{\left( 2,\alpha
\right) }\\
&=&\varepsilon \left( a_{\left( 1,1\right) }\right) a_{\left(
2,\alpha \right) }\varepsilon \left( b_{\left( 1,1\right) }\right)
d_{\alpha }b_{\left( 2,\alpha \right) }\\
&=&ad_{\alpha }b.
\end{eqnarray*}

\item Let $a\in A_{\alpha \beta }$
\begin{eqnarray*}
\Delta _{\alpha ,\beta }^{l}d_{\alpha \beta }\left( a\right)
&=&\Delta _{\alpha ,\beta }^{l}\left( d_{\alpha \beta }\left(
a\right) \right)\\
&=&\Delta _{\alpha ,\beta }\left( 1_{\alpha \beta }\right) \left(
id\otimes d_{\beta }\right)\Delta _{\alpha ,\beta
}\left( a\right)\\
&=&\left( 1_{\alpha }\otimes 1_{\beta }\right) \left( id\otimes
d_{\beta }\right)\Delta _{\alpha ,\beta }\left(
a\right)\\
&=&\left( id\otimes d_{\beta }\right)\Delta _{\alpha ,\beta
}\left( a\right)
\end{eqnarray*}
\end{enumerate}
\end{proof}
\begin{defn} Let $\Gamma =\left( \left\{ \Gamma_{\alpha}\right\} _{\alpha \in \pi },d\right) $ be a $\pi -$graded first order differential calculus over $A$.We say that $
\Gamma =\left( \left\{ \Gamma_{\alpha} \right\} _{\alpha \in \pi
},d\right) $ is right covariant if for any $\alpha ,\beta \in \pi
$ $$\sum_{k=1}^{n}a_{k}d_{\alpha \beta }b_{k}\Longrightarrow
\sum_{k=1}^{n}\Delta _{\alpha ,\beta }\left( a_{k}\right) \left(
d_{\alpha }\otimes id\right)\Delta _{\alpha ,\beta }\left(
b_{k}\right) =0\eqno{2.14}$$ \end{defn} \noindent We say that
$\Gamma =\left( \left\{ \Gamma_{\alpha} \right\} _{\alpha \in \pi
},d\right) $ is bicovariant if it is left and right covariant.
\begin{prop} Let $\Gamma =\left( \left\{ \Gamma_{\alpha} \right\}
_{\alpha \in \pi },d\right) $ be a right covariant $\pi -$graded
first order differential calculus over $A.$
 Then there exists a family of linear mappings $$\Delta ^{r}=\left\{ \Delta _{\alpha ,\beta }^{r}:\Gamma _{\alpha \beta }\longrightarrow \Gamma _{\alpha }\otimes A_{\beta }\right\} \eqno{2.15}$$
such that
\begin{enumerate}
\item For any $a\in A_{\alpha \beta },\rho \in \Gamma _{\alpha \beta }$ $$ \begin{array}{c} \Delta _{\alpha ,\beta }^{r}\left( a\rho \right) =\Delta _{\alpha ,\beta }\left( a\right) \Delta _{\alpha ,\beta }^{r}\left( \rho \right) \\
\Delta _{\alpha ,\beta }^{r}\left( \rho a\right) =\Delta _{\alpha ,\beta }^{r}\left( \rho \right) \Delta _{\alpha ,\beta }\left( a\right) \\
\end{array} \eqno{2.16} $$
\item for any $\alpha,\beta,\gamma\in\pi$ $$ (id\otimes
\Delta_{\beta,\gamma})\Delta^{r}_{\alpha,\beta\gamma}=(\Delta^{r}_{\alpha,\beta}
\otimes id)\Delta^{r}_{\alpha\beta,\gamma} \eqno{2.17}$$ \item For
any $\rho \in \Gamma _{\alpha }$ $$\left( id\otimes \varepsilon
\right) \Delta _{\alpha ,1}^{r}\left( \rho \right) =\rho
\eqno{2.18}$$ \item for any $\alpha,\beta,\gamma\in\pi$ $$
\Delta^{r}_{\alpha,\beta}d_{\alpha,\beta}=(d_{\alpha}\otimes
id)\Delta_{\alpha,\beta} $$ \end{enumerate}
\end{prop}
\begin{proof} Similar to that of proposition 2.2 , where for
any $\alpha,\beta,\gamma\in\pi$ $a_{k},b_{k}\in A_{\alpha \beta
}.$ $$\Delta _{\alpha ,\beta }^{r}\left( \sum_{k}a_{k}d_{\alpha
\beta }b_{k}\right) =\sum_{k}\Delta _{\alpha ,\beta }\left(
a_{k}\right) \left( d_{\alpha }\otimes id\right) \Delta _{\alpha
,\beta }\left( b_{k}\right) \eqno{2.19}$$
\end{proof}
\begin{prop} Let $\Gamma =\left( \left\{ \Gamma_{\alpha} \right\}
_{\alpha \in \pi },d\right) $ be a bicovariant $\pi -$graded first
order differential calculus over$A,$ $
 \Delta ^{l},\Delta ^{r}$ be the families of linear mappings introduced by
  proposition 2.2,2.3.Then we have. $$\left( id\otimes
  \Delta _{\beta ,\gamma }^{r}\right) \Delta _{\alpha ,\beta \gamma }^{l}
  \left( ad_{\alpha \beta \gamma }b\right) = \left( \Delta _{\alpha ,
  \beta }^{l}\otimes id\right) \Delta _{\alpha \beta ,\gamma }
  ^{r}\left( ad_{\alpha \beta \gamma }b\right)\eqno{2.20}$$
\end{prop}
\begin{proof} Let $a,b\in A_{\alpha \beta \gamma }$
\begin{eqnarray*}
\left( id\otimes \Delta _{\beta ,\gamma }^{r}\right) \Delta
_{\alpha ,\beta \gamma }^{l}\left( ad_{\alpha \beta \gamma
}b\right)
 &=&\left( id\otimes \Delta _{\beta ,\gamma }^{r}\right)
  \left( \Delta _{\alpha ,\beta \gamma }\left( a\right)
   \left( id\otimes d_{\beta \gamma }\right)
   \Delta _{\alpha ,\beta \gamma }\left( b\right) \right)\\
   &=&\left( id\otimes \Delta _{\beta ,\gamma }\right)
    \Delta _{\alpha ,\beta \gamma }\left( a\right)
    \left( id\otimes d_{\beta }\otimes id
    \right) \left( id\otimes \Delta _{\beta ,\gamma }\right)
     \Delta _{\alpha ,\beta \gamma }\left( b\right)
     \end{eqnarray*}
On the other hand
\begin{eqnarray*} \left( \Delta _{\alpha ,\beta
}^{l}\otimes id\right) \Delta _{\alpha \beta ,\gamma }^{r}\left (
ad_{\alpha \beta \gamma }b\right) &=& \left( \Delta_{\alpha ,
\beta }^{l}\otimes id\right) \left( \Delta _{\alpha \beta ,\gamma
}\left( a\right) \left( id\otimes d_{\gamma }\right) \Delta
_{\alpha \beta ,\gamma
}\left( b\right) \right)\\
&=&\left( \Delta _{\alpha ,\beta}\otimes id\right) \Delta _{\alpha
\beta ,\gamma }\left( a\right) \left( id\otimes d_{\beta }\otimes
id\right) \left( \Delta _{\alpha ,\beta }\otimes id\right) \Delta
_{\alpha ,\beta \gamma }\left( b\right)
\end{eqnarray*}
Using the coassociativity property we find that equation 2.20
holds.
 \end{proof}
\section{$\pi-$graded Bicovariant bimodules}
Throughout this section let $A=$ $\left\{
A_{\alpha}\right\}_{\alpha \in \pi }$ be a hopf group coalgebra
%==============================================================================================================
\begin{defn}
let $\Gamma =\left\{ \Gamma _{\alpha}\right\}_{\alpha \in \pi }$
be a $\pi -$graded bimodule over A , $\Delta ^{l}=\left\{ \Delta
_{\alpha ,\beta }^{l}:\Gamma _{\alpha \beta }\longrightarrow
A_{\alpha }\otimes \Gamma _{\beta }\right\} _{\alpha ,\beta \in
\pi }$ be a family of linear maps .We say that $\Gamma =\left(
\left\{ \Gamma _{\alpha }\right\}_{\alpha \in \pi } ,\Delta
^{l}\right) $ is a left covariant $\pi $-graded bimodule over A if
\begin{enumerate}
    \item For any $a\in A_{\alpha \beta} \;\;,\rho \in \Gamma _{\alpha \beta} \;\; \alpha ,\beta \in \pi $

$$\Delta _{\alpha ,\beta }^{l}\left( a\rho \right) =\Delta _{\alpha
,\beta }\left( a\right) \Delta _{\alpha ,\beta }^{l}\left( \rho
\right) \eqno(3.1)$$

$$ \Delta _{\alpha ,\beta }^{l}\left( \rho
a\right) =\Delta _{\alpha ,\beta }^{l}\left( \rho \right) \Delta
_{\alpha ,\beta }\left( a\right) \eqno(3.2)$$

    \item For all  $\alpha ,\beta ,\gamma \in \pi .$
    $$
(\Delta_{\alpha,\beta}\otimes id)\Delta_{\alpha\beta,\gamma}^{l}=
(id\otimes
\Delta_{\beta,\gamma}^{l})\Delta_{\alpha,\beta\gamma}^{l}
\eqno(3.3)
$$
    \item For any $\rho \in \Gamma _{\alpha }$,$\alpha \in \pi $

$$ \left( \varepsilon \otimes id\right) \Delta _{1,\alpha }^{l}\left( \rho \right) =\rho \eqno
(3.4)$$
\end{enumerate}
\end{defn}

%===============================================================================================================
\begin{defn}
Let $\Gamma =\left\{ \Gamma _{\alpha }\right\} _{\alpha \in \pi
}$be a $\pi -$graded bimodule over A ,
 $\Delta ^{r}=\left\{ \Delta _{\alpha ,\beta }^{r}:\Gamma _{\alpha \beta }\longrightarrow \Gamma _{\alpha }\otimes
 A_{\beta }\right\}$ be a family of linear maps .We say that $\Gamma =\left( \left\{ \Gamma
_{\alpha }\right\} _{\alpha \in \pi },\Delta ^{r}\right) $is a
right covariant $\pi -$graded bimodule over A if
\begin{enumerate}
    \item For any $a\in A_{\alpha \beta },\rho \in \Gamma _{\alpha \beta }$

$$ \Delta _{\alpha ,\beta }^{r}\left( a\rho
\right) =\Delta _{\alpha ,\beta }\left( a\right) \Delta _{\alpha
,\beta }^{r}\left( \rho \right) \eqno(3.5)$$

$$ \Delta _{\alpha ,\beta }^{r}\left( \rho a\right) =\Delta
_{\alpha ,\beta }^{r}\left( \rho \right) \Delta _{\alpha ,\beta
}\left( a\right)\eqno(3.6)$$
    \item For $\alpha ,\beta ,\gamma \in \pi .$

  $$
(\Delta_{\alpha,\beta}\otimes id)\Delta_{\alpha\beta,\gamma}^{r}=
(id\otimes
\Delta_{\beta,\gamma}^{r})\Delta_{\alpha,\beta\gamma}^{r}
\eqno(3.7)
$$

    \item For any $\rho \in \Gamma _{\alpha }$,$\alpha \in \pi $

$$\left( id\otimes \varepsilon
\right) \Delta _{\alpha ,1}^{r}\left( \rho \right) =\rho \eqno
(3.8)$$
\end{enumerate}

\end{defn}
%=============================================================================================
\begin{defn}
let $\Gamma =\left\{ \Gamma _{\alpha }\right\}_{\alpha \in \pi } $
be a $\pi -$graded bimodule over A , $\Delta ^{l}=\left\{ \Delta
_{\alpha ,\beta }^{l}:\Gamma _{\alpha \beta }\longrightarrow
A_{\alpha }\otimes \Gamma _{\beta }\right\} _{\alpha ,\beta \in
\pi }$, and $\Delta ^{r}=\left\{ \Delta _{\alpha ,\beta
}^{r}:\Gamma _{\alpha \beta }\longrightarrow \Gamma _{\alpha
}\otimes A_{\beta }\right\} $be two families of linear maps .We
say that $\Gamma =\left( \left\{ \Gamma _{\alpha }\right\}_{\alpha
\in \pi } ,\Delta ^{l},\Delta ^{r}\right) $ is a bicovariant $\pi
$-graded bimodule over A if
\begin{enumerate}
    \item $\Gamma =\left( \left\{ \Gamma _\alpha
    \right\}_{\alpha \in \pi }
,\Delta ^{l}\right) $ is a left covariant $\pi $-graded bimodule
over A.
    \item $\Gamma =\left( \left\{ \Gamma _\alpha \right\}_{\alpha \in \pi }
,\Delta ^{r}\right) $ is a right covariant $\pi $-graded bimodule
over A.
    \item For all $\alpha ,\beta ,\gamma \in \pi .$
$$
(\Delta_{\alpha,\beta}^{l}\otimes
id)\Delta_{\alpha\beta,\gamma}^{r}= (id\otimes
\Delta_{\beta,\gamma}^{r})\Delta_{\alpha,\beta\gamma}^{l}
\eqno(3.9)
$$
\end{enumerate}
\end{defn}
%==================================================================================================
\begin{defn}
Let $\Gamma =\left( \left\{ \Gamma _\alpha
    \right\}_{\alpha \in \pi },\Delta ^{l}\right) $ be a left covariant $\pi $-graded bimodule
over A .For any $\alpha \in \pi$ an element $\rho \in \Gamma
_{\alpha }$ is said to be left invariant if

$$ \Delta _{1,\alpha }^{l}\left( \rho \right) =1_{1}\otimes \rho
\eqno(3.10)$$
\end{defn}

\medskip

Denote by $_{inv}\Gamma =\left\{ _{inv}\Gamma _{\alpha }\right\}
_{\alpha \in \pi }$ the set of all left invariant elements of
$\Gamma $. Clearly , $ _{inv}\Gamma _{\alpha }$ is a linear
subspace of $\Gamma _{\alpha }$ for each $\alpha \in \pi .$

\medskip

%==========================================================================================
\begin{lem}
Let $\Gamma =\left( \left\{ \Gamma _{\alpha }\right\} _{\alpha \in
\pi
},\Delta ^{l}\right) $ be a left covariant $\pi $-graded bimodule over $%
A,_{inv}\Gamma =\left\{ _{inv}\Gamma _{\alpha }\right\} _{\alpha
\in \pi }$ be the linear subspace of all left invariant elements
of $\Gamma$. Then there exists a family

$$P=\left\{ P_{\alpha }:\Gamma
_{1}\longrightarrow \Gamma _{\alpha }\right\} _{\alpha \in \pi
}\eqno(3.11)$$

of mappings such that

$$ P_{\alpha }\left( b\rho \right)
=\varepsilon \left( b\right) P_{\alpha }\left( \rho \right) \eqno
(3.12)$$

for any b$\in A_{1},\rho \in \Gamma _{1},\alpha \in \pi .$

Moreover , for any $\rho \in \Gamma _{\alpha }$ $,\alpha \in \pi $
we have
$$ \rho =\sum_{k}a_{k}P_{\alpha }\left( \rho _{k}\right)
\eqno(3.13) $$

where $a_{k}$ , $\rho _{k}$ are elements of $A_{\alpha }$ ,
$\Gamma _{1}$ respectively such that

$$ \Delta _{\alpha ,1}^{l}\left( \rho \right) =\sum_{k}a_{k}\otimes \rho _{k}\eqno(3.14)$$

And equation 3.13 can be uniquely written in this form.
\end{lem}
%==================================================================================================================
\begin{proof}

For any $\alpha \in \pi $ , $\rho \in \Gamma _{1}$ set

$$P_{\alpha }\left( \rho \right) =\sum_{k}S_{\alpha ^{-1}}\left(
a_{k}\right) \rho _{k} \eqno(3.15)$$

where

$$\Delta _{\alpha ^{-1},\alpha }^{l}\left(
\rho \right) =\sum_{k=1}^{n}a_{k}\otimes \rho _{k}$$

Recall that for any $\alpha ,\beta \in \pi $ , $a\in A_{\beta
^{-1}}$ where  $ \Delta _{\beta ^{-1}\alpha ^{-1},\alpha }\left(
a\right) = a_{\left( 1,\beta ^{-1}\alpha ^{-1}\right) }\otimes
a_{\left( 2,\alpha \right) }$we have

%==================================================================================================================
\begin{eqnarray*}
\Delta _{\alpha ,\beta }\left( S_{\beta ^{-1}\alpha ^{-1}}\left(
a_{\left( 1,\beta ^{-1}\alpha ^{-1}\right) }\right) \right) \left(
a_{\left( 2,\alpha \right) }\otimes 1_{\beta }\right)&=&\sigma
_{A_{\beta },A_{\alpha }}\left( S_{\beta ^{-1}}\left( a_{\left(
1,\beta ^{-1}\right) }\right) \otimes S_{\alpha ^{-1}}\left(
a_{\left( 2,\alpha ^{-1}\right) }\right) \right)
\left( a_{\left( 3,\alpha \right) }\otimes 1_{\beta }\right)\\
&=&S_{\alpha ^{-1}}\left( a_{\left( 2,\alpha ^{-1}\right) }\right)
a_{\left( 3,\alpha \right) }\otimes S_{\beta ^{-1}}\left(
a_{\left( 1,\beta ^{-1}\right) }\right)\\
&=&\varepsilon \left( a_{\left( 2,1\right) }\right) 1_{\alpha
}\otimes S_{\beta ^{-1}}\left( a_{\left( 1,\beta ^{-1}\right)
}\right)\\
&=&1_{\alpha }\otimes \varepsilon \left( a_{\left( 2,1\right)
}\right) S_{\beta ^{-1}}\left( a_{\left( 1,\beta ^{-1}\right)
}\right)\\
&=& 1_{\alpha }\otimes S_{\beta ^{-1}}\left( a_{\left( 1,\beta
^{-1}\right) }\varepsilon \left( a_{\left( 2,1\right) }\right)
\right)\\
&=& 1_{\alpha }\otimes S_{\beta ^{-1}}\left( a \right)
\end{eqnarray*}
then we have
$$ \Delta _{\alpha ,\beta }\left( S_{\beta ^{-1}\alpha ^{-1}}\left(
a_{\left( 1,\beta ^{-1}\alpha ^{-1}\right) }\right) \right) \left(
a_{\left( 2,\alpha \right) }\otimes 1_{\beta }\right)=1_{\alpha
}\otimes S_{\beta ^{-1}}\left( a\right) \eqno(3.16)$$
%==================================================================================================================
For any $\rho \in \Gamma _{1\textbf{ , }}\alpha \in \pi $ set
\begin{eqnarray*}
\Delta _{\alpha ^{-1},\alpha }^{l}\left( \rho \right)
&=&\sum_{k}a_{k}\otimes \rho _{k}\\
\Delta _{1,\alpha }^{l}\left( \rho _{k}\right)
&=&\sum_{l}b_{kl}\otimes \rho _{kl}\\
\Delta _{\alpha ^{-1},1}\left( a_{k}\right)
&=&\sum_{m}c_{km}\otimes d_{km}
\end{eqnarray*}

Using equation 3.3 we have

$$ \sum_{k,l}a_{k}\otimes b_{kl}\otimes \rho
_{kl}=\sum_{k,m}c_{km}\otimes d_{km}\otimes \rho _{k}
\eqno(3.17)$$

We compute
%==================================================================================================================
\begin{eqnarray*}
\Delta _{1,\alpha }^{l}\left( P_{\alpha }\left( \rho \right)
\right) &=&\Delta _{1,\alpha }^{l}\left( \sum_{k}S_{\alpha
^{-1}}\left( a_{k}\right) \rho _{k}\right) \\
&=&\sum_{k}\Delta _{1,\alpha }\left( S_{\alpha ^{-1}}\left(
a_{k}\right) \right) \Delta _{1,\alpha }^{l}\left( \rho
_{k}\right) \\
&=&\sum_{k,l}\Delta _{1,\alpha }\left( S_{\alpha ^{-1}}\left(
a_{k}\right) \right) \left( b_{kl}\otimes \rho _{kl}\right) \\
&=&\sum_{k,l}\Delta _{1,\alpha }\left( S_{\alpha ^{-1}}\left(
c_{km}\right) \right) \left( d_{km}\otimes \rho _{k}\right) \\
&=&\sum_{k,l}\Delta _{1,\alpha }\left( \otimes S_{\alpha
^{-1}}\left( c_{km}\right) \right) \left( d_{km}\otimes 1_{\alpha
}\right) \left( 1_{1}\otimes \rho _{k}\right) \\
&=&\sum_{k}\left( 1_{1}\otimes S_{\alpha ^{-1}}\left( a_{k}\right)
\right) \left( 1_{1}\otimes \rho _{k}\right) \\
\left( \textrm{using 3.16 for }\alpha =1,\beta =\alpha \right)
\\
&=&1_{1}\otimes \sum_{k}S_{\alpha ^{-1}}\left( a_{k}\right) \rho
_{k}\\
&=&1_{1}\otimes P_{\alpha }\left( \rho \right) \\
\end{eqnarray*}
This shows that $P_{\alpha }\left( \rho \right) $ is left
invariant element in $\Gamma _{\alpha }$ for each $\alpha \in \pi
.$

To prove 3.12 , let $b\in A_{1}$ , $\rho \in \Gamma _{1}$ , set
\begin{eqnarray*}
\Delta _{\alpha ^{-1},\alpha }\left(
b\right)&=&\sum_{k}b_{k}\otimes d_{k} \\
\Delta _{\alpha ^{-1},\alpha }^{l}\left( \rho \right)
&=&\sum_{l}c_{l}\otimes \rho _{l}\\
\Delta _{\alpha ^{-1},\alpha }^{l}\left( b\rho \right) &=&\Delta
_{\alpha ^{-1},\alpha }\left( b\right) \Delta _{\alpha
^{-1},\alpha }^{l}\left( \rho \right) \\
&=& \sum_{k,l}b_{k}c_{l}\otimes d_{k}\rho _{l}
\end{eqnarray*}
Then
%=====================================================================================================
\begin{eqnarray*}
P_{\alpha }\left( b\rho \right) &=&\sum_{k,l}S_{\alpha
^{-1}}\left(
b_{k}c_{l}\right) d_{k}\rho _{l}\\
&=&\sum_{k,l}S_{\alpha ^{-1}}\left( c_{l}\right)
S_{\alpha ^{-1}}\left( b_{k}\right) d_{k}\rho _{l}\\
&=&\sum_{l}S_{\alpha ^{-1}}\left( c_{l}\right)
\varepsilon \left( b\right) \rho _{l}\\
&=&\varepsilon \left( b\right) \sum_{l}S_{\alpha
^{-1}}\left( c_{l}\right) \rho _{l}\\
&=&\varepsilon \left( b\right) P_{\alpha }\left( \rho \right)
\end{eqnarray*}
To prove 3.13. Let $\alpha \in \pi ,\rho \in \Gamma _{\alpha
}.$Set
\begin{eqnarray*}
\Delta _{1,\alpha }^{l}\left( \rho \right) &=&\sum_{m}d_{m}\otimes
\varrho _{m} \\
\Delta _{\alpha ^{-1},\alpha }^{l}\left( \rho _{k}\right)
&=&\sum_{n}b_{kn}\otimes \rho _{kn}\\
 \Delta _{\alpha ,\alpha
^{-1}}\left( d_{m}\right) &=&\sum_{l}d_{ml}\otimes c_{ml}\\
\end{eqnarray*}
where
$$\Delta _{\alpha ,1}^{l}\left( \rho \right) =\sum_{k}a_{k}\otimes
\rho _{k} \eqno(3.18) $$ using equation 3.3 we have

$$\left( \Delta _{\alpha ,\alpha ^{-1}}\otimes id\right) \Delta _{1,\alpha }
^{l}=\left( id\otimes \Delta _{\alpha ^{-1},\alpha }^{l}\right)
\Delta _{\alpha ,1}^{l}$$

i.e.
$$\sum_{m,l}d_{ml}\otimes c_{ml}\otimes \varrho
_{m}=\sum_{k,n}a_{k}\otimes b_{kn}\otimes \rho _{kn}\eqno(3.19) $$
Then using equation 3.4 we have
\begin{eqnarray*}
\rho &=&\left( \varepsilon \otimes id\right) \Delta _{1,\alpha
}^{l}\left( \rho \right)
\\
&=&\sum_{m}\varepsilon \left(
d_{m}\right) \varrho _{m}\\
&=&\sum_{m,l}d_{ml}S_{\alpha
^{-1}}\left( c_{ml}\right) \varrho _{m}\\
&=&\sum_{k,n}a_{k}S_{\alpha ^{-1}}\left(
b_{kn}\right) \rho _{kn}\\
&=&\sum_{m}a_{k}P_{\alpha }\left( \rho _{k}\right)
\end{eqnarray*}
Finally ,to prove the uniqueness of expression 3.13 let
$P^{^{\prime }}=\left\{ P_{\alpha }^{^{\prime }}:\Gamma
_{1}\longrightarrow _{inv}\Gamma _{\alpha }\right\} $ be another
family of mappings satisfying that for $\rho\in\Gamma_\alpha$
 $$ \rho =\sum_{k}a_{k}P_{\alpha }^{\prime
}\left( \rho _{k}\right)\eqno(3.20)$$

where $a_{k}$ , $\rho _{k}$ are elements of $A_{\alpha }$ ,
$\Gamma _{1}$ respectively such that

$$\Delta _{\alpha ,1}^{l}\left( \rho
\right) =\sum_{k}a_{k}\otimes \rho _{k}$$

Let $\rho \in \Gamma _{\alpha }$ such that $\Delta _{\alpha
,1}^{l}\left( \rho \right) =\sum_{k}a_{k}\otimes \rho _{k}$ .Then
using 3.13

$$ \rho =\sum_{k}a_{k}P_{\alpha }\left( \rho _{k}\right)$$
But using 3.20 we have
$$ \rho =\sum_{k}a_{k}P_{\alpha }^{^{\prime }}\left( \rho
_{k}\right) $$ Subtracting the above two equations we obtain
$$ 0=\sum_{k}a_{k}\left( P_{\alpha }\left( \rho _{k}\right)
-P_{\alpha }^{^{\prime }}\left( \rho _{k}\right) \right) $$
Assuming that all $a_{k}^{,}s$ all linearly independent we get
$$ P_{\alpha }\left( \rho _{k}\right) =P_{\alpha
}^{^{\prime }}\left( \rho _{k}\right)
\;\;\;\;\;\;\;\;k=1,2,...,n.$$ which proves the uniqueness of
expression 3.13.

\begin{lem}
Let $\Gamma =\left( \left\{ \Gamma _{\alpha }\right\} _{\alpha \in
\pi },\Delta ^{l}\right) $ be a left covariant $\pi $-graded
bimodule over A . Then , for any $\alpha ,\beta \in \pi ,\rho \in
_{inv}\Gamma _{\alpha \beta } $ we have

$$\Delta _{\alpha ,\beta
}^{l}\left( \rho \right) =1_{\alpha }\otimes \varrho \eqno(3.21)
$$ where $\varrho \in _{inv}\Gamma _{\beta }.$
\begin{proof}
Let $\alpha ,\beta \in \pi ,\rho \in _{inv}\Gamma _{\alpha \beta
}$ ,then using lemma 3.1 , and since the mappings $P_{\alpha }$
are onto for each $\alpha \in \pi $ then there exists an element
$\xi \in \Gamma_{1}$ such that
$$\rho =P_{\alpha \beta }\left( \xi \right)\eqno( 3.22)$$
Set
\begin{eqnarray*}
\Delta _{\beta ^{-1}\alpha ^{-1},\alpha \beta }^{l}\left( \xi
\right) &=&\sum_{k}a_{k}\otimes
\xi _{k}\\
\Delta _{\alpha ,\beta }^{l}\left(
\xi _{k}\right) &=&\sum_{l}c_{kl}\otimes \xi _{kl}\\
\Delta _{\beta ^{-1},\beta }^{l}\left( \xi
\right)&=&\sum_{m}b_{m}\otimes \rho _{m}                 \\
\end{eqnarray*}
and
 $$ \Delta _{\beta ^{-1}\alpha ^{-1},\alpha }^{l}\left(
b_{m}\right) =\sum_{n}b_{mn}\otimes d_{mn}\eqno( 3.23)$$

 Using equation 3.3
$$\sum_{k,l}a_{k}\otimes c_{kl}\otimes \xi
_{kl}=\sum_{m,n}b_{mn}\otimes d_{mn}\otimes \rho _{m}
\eqno(3.24)$$ Applying $\Delta _{\alpha ,\beta }^{l}$ to both
sides of 3.22, using 3.24 and 3.16 , we get
\begin{eqnarray*}
\Delta _{\alpha ,\beta }^{l}\left( \rho \right) &=&\Delta _{\alpha
,\beta }^{l}\left( P_{\alpha \beta }\left( \xi
\right) \right) \\
&=&\sum_{k}\Delta _{\alpha ,\beta }^{l}\left( S_{\beta ^{-1}\alpha
^{-1}}\left( a_{k}\right) \xi
_{k}\right)\\
&=&\sum_{k}\Delta _{\alpha ,\beta }\left( S_{\beta ^{-1}\alpha
^{-1}}\left( a_{k}\right) \right)
\Delta _{\alpha ,\beta }^{l}\left( \xi _{k}\right)\\
&=&\sum_{k,l}\Delta _{\alpha ,\beta }\left( S_{\beta ^{-1}\alpha
^{-1}}\left( a_{k}\right) \right)
\left( c_{kl}\otimes \xi _{kl}\right) \\
&=&\sum_{m,n}\Delta _{\alpha ,\beta }\left( S_{\beta ^{-1}\alpha
^{-1}}\left( b_{mn}\right) \right)
\left( d_{mn}\otimes \rho _{m}\right) \\
&=&\sum_{m,n}\Delta _{\alpha ,\beta }\left( S_{\beta ^{-1}\alpha
^{-1}}\left( b_{mn}\right) \right) \left( d_{mn}\otimes 1_{\beta
}\right) \left( 1_{\alpha }\otimes
\rho _{m}\right) \\
&=&\sum_{m}\left( 1_{\alpha }\otimes S_{\beta ^{-1}}\left(
b_{m}\right) \right) \left( 1_{\alpha
}\otimes \rho _{m}\right) \\
&=&1_{\alpha }\otimes \sum_{m}S_{\beta
^{-1}}\left( b_{m}\right) \rho _{m}\\
&=&1_{\alpha }\otimes P_{\beta }\left( \xi \right)
\end{eqnarray*}
But from lemma 3.1 $P_{\beta } ( \xi ) \in _{inv} \Gamma _{\beta
}$ and hence the lemma is proved.
\end{proof}
\end{lem}

Let $A= ( \{ A_{\alpha } \}_{\alpha\in\pi},\Delta ,\varepsilon,S
)$ be a hopf $\pi -$coalgebra. Throughtout the next dealing we
will consider that $A$ is endowed with a family of linear maps
$\Psi = \{ \Psi _{\alpha }: A_{\alpha } \longrightarrow A_{1} \}$
of k-linear maps such that for each $\alpha \in \pi $ ,$\Psi
_{\alpha }$ is an algebra map.For each $\alpha \in \pi $ ,define
the map $E_{\alpha }$ to be the composition
$$ A_{\alpha }\longrightarrow
A_{1}\longrightarrow k$$ i.e.
$$ E_{\alpha }=\varepsilon \Psi _{\alpha }\eqno(3.25)$$

Clearly, for each $\alpha \in \pi $ $E_{\alpha }$ is an algebra
map for let $a,b\in A_{\alpha }$ .Then

\begin{eqnarray*}
E_{\alpha }\left( ab\right) &=&\varepsilon \left( \Psi
_{\alpha }\left( ab\right) \right) \\
&=&\varepsilon \left( \Psi _{\alpha }\left(
a\right) \textrm{ }\Psi _{\alpha }\left( b\right) \right) \\
&=&\varepsilon \left( \Psi _{\alpha }\left( a\right) \textrm{
}\right) \varepsilon \left( \textrm{ }\Psi
_{\alpha }\left( b\right) \right) \\
&=&E_{\alpha }\left( a\right)E_{\alpha
}\left( b\right)\\
E_{\alpha }\left( 1_{\alpha }\right) &=&\varepsilon \left(
\Psi _{\alpha }\left( 1_{\alpha }\right) \right) \\
&=&\varepsilon \left( 1_{1}\right) \\
&=&1_{k}
\end{eqnarray*}
Moreover, $E_{\alpha }$ is linear being the composition of two
linear maps.

\begin{thm}
Let $\Gamma =\left( \left\{ \Gamma _{\alpha }\right\} _{\alpha \in
\pi},\Delta ^{l}\right) $ be a $\pi -$graded left covariant bimodule over A,$%
\left\{ \omega _{i}^{\alpha }\right\} _{\alpha \in \pi }$be a
basis of $_{inv}\Gamma _{\alpha },$ of all left invariant elements
of $\Gamma _{\alpha }$ for each $\alpha \in \pi $ .Then
\begin{enumerate}
\item For any $\alpha \in \pi $ ,any element $\rho \in $ $\Gamma
_{\alpha }$ is of the form
$$\rho = \sum_{i}a_{i}\omega _{i} \eqno(3.26)$$
where $a_{i}$ $^{,}s\in A_{\alpha }$ are uniquely determined ,
$\omega _{i}^{,}s\in _{inv}\Gamma _{\alpha }$ ,for any $\alpha \in
\pi $ . \item For any $\alpha \in \pi $ ,any element $\rho \in $
$\Gamma _{\alpha }$ is of the form
$$\rho =\sum_{i}\omega
_{i}b_{i}\eqno(3.27)$$ where $b_{i}$ $^{,}s\in A_{\alpha }$ are
uniquely determined , $\omega _{i}^{,}s\in _{inv}\Gamma _{\alpha
}$ ,for any $\alpha \in \pi $ . \item There exists linear
functionals $f_{ij}$ $\in A^{^{\prime }}=\oplus _{\alpha \in \pi
}A_{\alpha }^{^{\prime }}$ such that for any $\alpha \in \pi $
$$\omega_{i}b=\sum_{j}\left( f_{ij}*b\right) \omega _{j}\eqno(3.28)$$
$$a\omega_{i}=\sum_{j}\omega _{j}\left( \left( f_{ij}\circ
S_{1}^{-1}\right) *a\right)\eqno(3.29)$$ where $a,b\in A_{\alpha
},\omega _{i}^{,}s,\omega _{j}^{,}s\in _{inv}\Gamma _{\alpha }$
.These functionals are uniquely determined by $3.28.\ $They
satisfy the following relations
$$f_{ij}(ab)=\sum_{k}f_{ik}(a)f_{kj}(b)\eqno(3.30)$$
for any $i,j\in I$ , $a,b\in A_{\alpha }$ .Moreover
$$f_{ij}(1_{\alpha})=\delta _{ij}\eqno(3.31)$$
\end{enumerate}
\end{thm}

\begin{rem}
Any functional $f_{ij}\in A^{^{\prime }}=\oplus _{\alpha \in \pi
}A_{\alpha }^{^{\prime }}$ is of the form $f_{ij}=\sum_{\alpha }$
$f_{ij}^{\alpha }$ where
$$f_{ij}^{\alpha }\left(a\right) =0 \;\;\; if \; a\notin A_{\alpha }$$
\end{rem}

\begin{proof}
To prove 1:
 For any $\alpha \in \pi $ let $\rho \in \Gamma
_{\alpha }$. Using 3.13 we have that $\rho =\sum_{i}a_{i}\omega
_{i}$ , with $\omega _{i}^{,}s\in _{inv}\Gamma _{\alpha }$ . To
prove uniqueness assume that $\rho =\sum_{i}a_{i}\omega _{i}$.
Then , using $3.1 , 3.21$
\begin{eqnarray*}
\Delta _{\alpha ,1}^{l}\left( \rho \right) &=&\Delta _{\alpha
,1}^{l}\left( \sum_{i}a_{i}\omega _{i}\right)
\\
&=&\sum_{i}\Delta _{\alpha,1}\left( a_{i}\right) \Delta _{\alpha
,1}^{l}\left( \omega _{i}\right)
\\
&=&\sum_{i}\left( a_{i\left(1,\alpha \right) }\otimes a_{i\left(
2,1\right) }\right) \left( 1_{\alpha }\otimes \xi _{i}\right) ,\xi
_{i}\in _{inv}\Gamma _{1}
\\
&=&\sum_{i}a_{i\left( 1,\alpha\right) }\otimes a_{i\left(
2,1\right) }\xi _{i}
\end{eqnarray*}
Applying $\left( id\otimes P_{1}\right) $ to both sides of the
above equation ,we get
\begin{eqnarray*}
\left( id\otimes P_{1}\right) \Delta _{\alpha ,1}^{l}\left( \rho
\right) &=&\left( id\otimes P_{1}\right) \sum_{i}a_{i\left(
1,\alpha \right) }\otimes a_{i\left( 2,1\right) }\xi _{i}
\\
&=&\sum_{i}a_{i\left(1,\alpha \right) }\otimes P_{1}\left(
a_{i\left( 2,1\right) }\xi _{i}\right)
\\
&=&\sum_{i}a_{i\left( 1,\alpha \right) }\otimes \varepsilon \left(
a_{i\left( 2,1\right) }\right) P_{1}\left( \xi _{i}\right)
\\
&=&\sum_{i}a_{i\left( 1,\alpha \right) }\varepsilon \left(
a_{i\left( 2,1\right) }\right) \otimes P_{1}\left( \xi _{i}\right)
\\
&=&\sum_{i}a_{i}\otimes P_{1}\left( \xi _{i}\right)
\\
&=&\sum_{i}a_{i}\otimes\xi _{i}
\end{eqnarray*}
since $P_{1}\left( \xi _{i}\right) =\xi _{i}$ for any $\xi _{i}\in
_{inv}\Gamma _{1}$ .Since $\omega _{i}^{,}s$ ,$i\in I$ are
linearly independent , then by linearty of $\Delta _{\alpha
,1}^{l}$ , $\xi _{i}^{,}s$ are also linearly independent and so
the coefficients $a_{i}^{,}s$ are uniquely determined, and this
proves the uniqueness of the decomposition $3.26.$ To prove 3: For
any $\alpha \in \pi $ , let $b\in A_{\alpha },\omega _{j}\in
_{inv}\Gamma _{\alpha } , j\in I$ . Then $\omega _{j}b$ admits a
decomposition in the of the form 3.26 .Let $F_{ji}^{\alpha }\left(
b\right) $ be the coefficients preceding $\omega _{i}$ in the
decomposition 3.26 i.e.
$$\omega _{j}b=\sum_{i}F_{ji}^{\alpha }\left(
b\right) \omega _{i}\eqno(3.32)$$ Clearly, $F_{ji}^{\alpha }\left(
b\right) $ are linear mappings acting on $A_{\alpha }$. For any
$a,b\in $ $A_{\alpha }$ ,and any $j\in I$ we have
\begin{eqnarray*}
\sum_{i}F_{ji}^{\alpha }\left( ab\right) \omega_{i}&=&\omega
_{j}ab
\\
&=&\sum_{h}F_{jh}^{\alpha}\left( a\right) \omega _{h}b
\\
&=&\sum_{h,i}F_{jh}^{\alpha}\left( a\right) F_{hi}^{\alpha }\left(
b\right) \omega _{i}
\end{eqnarray*}
using the uniqueness of the decomposition 3.26 we have
$$F_{ji}^{\alpha }\left( ab\right)
=\sum_{h}F_{jh}^{\alpha }\left( a\right) F_{hi}^{\alpha }\left(
b\right) \eqno(3.33)$$ for all $i,j\in I,\alpha \in \pi ,a,b\in
A_{\alpha }.$ Let $f_{ji\textrm{ }}^{\alpha }$be linear
functionals defined on $A_{\alpha }$ introduced by the formula
$$f_{ji\textrm{ }}^{\alpha }\left( a\right)
=E_{\alpha }\left( F_{ji}^{\alpha }\left( a\right) \right)
=\varepsilon \left( \Psi _{\alpha }\left( F_{ji}^{\alpha }\left(
a\right) \right) \right) \eqno(3.34)$$ Define $f_{ji}\in
A^{^{\prime }}$ by
$$f_{ji}=\sum_{\alpha \in \pi }f_{ji\textrm{ }}^{\alpha }$$
where for any $\beta \in \pi ,a\in A_{\beta }$
$$f_{ji}\left( a\right)=\sum_{\alpha \in \pi } f_{ji\textrm{ }}^{\alpha }\left(a\right) =f_{ji}^{\beta }\left( a\right) \eqno(3.35)$$
Applying $E_{\alpha }$ to both sides of 3.33 and using 3.34 and
3.35 we have
$$f_{ji}\left( ab\right)=\sum_{h}f_{jh}\left( a\right) f_{hi}\left( b\right) $$
for any $a,b\in A_{\alpha }$, and hence 3.30 is proven. From 3.30
we get
$$f_{ji}m_{\alpha }\left( a\otimes b\right) =\sum_{h}\left( f_{jh}\otimes f_{hi}\right) \left(
a\otimes b\right) $$ i.e.
$$f_{ji}m_{\alpha }=\sum_{h}\left(f_{jh}\otimes f_{hi}\right) \eqno(3.36)$$
Inserting $b=1_{\alpha }$ in 3.32 we get

$$\omega _{j}=\sum_{i}F_{ji}^{\alpha
}\left( 1_{\alpha }\right) \omega _{i}$$ i.e.
$$F_{ji}^{\alpha }\left( 1_{\alpha
}\right) =\delta _{ji}1_{\alpha }$$ Applying $E_{\alpha }$ to both
sides of the above equation ,and summing over $\alpha $ we get
$$f_{ji}\left( 1_{\alpha }\right) =\delta _{ji}$$
and hence 3.31 is proven. To prove 3.28 Recall that from equation
3.32 for any $\alpha \in \pi,\omega _{j}\in_{inv}\Gamma _{\alpha
},$ $b\in A_{\alpha }$
$$\omega _{j}b=\sum_{i}F_{ji}^{\alpha}\left( b\right) \omega _{i}$$
Applying $\Delta _{\alpha ,1}^{l}$ to both sides of the above
equation we obtain
\begin{eqnarray*}
\Delta _{\alpha ,1}^{l}\left( \omega_{j}b\right) &=&\Delta
_{\alpha ,1}^{l}\left( \sum_{i}F_{ji}^{\alpha}\left( b\right)
\omega _{i}\right)
\\
\left( 1_{\alpha }\otimes \xi _{j}\right) \Delta _{\alpha
,1}\left( b\right) &=&\sum_{i}\Delta _{\alpha ,1}\left(
F_{ji}^{\alpha }\left( b\right) \right) \left( 1_{\alpha }\otimes
\xi _{i}\right)
\end{eqnarray*}
where $\xi _{j},\xi _{i}\in _{inv}\Gamma _{1}$ , $i,j\in I.$ On
the other hand using 3.32
$$\left( 1_{\alpha }\otimes \xi _{j}\right) \Delta _{\alpha
,1}\left( b\right) =\sum_{i}\left( id\otimes F_{ji}^{1}\right)
\Delta _{\alpha ,1}\left( b\right) \left( 1_{\alpha }\otimes \xi
_{i}\right) $$ Comparing the last two equations we get
$$\Delta _{\alpha ,1}\left(
F_{ji}^{\alpha }\left( b\right) \right) =\left( id\otimes
F_{ji}^{1}\right) \Delta _{\alpha ,1}\left( b\right) $$ Applying
$\left( id\otimes \varepsilon \right) $ to both sides of the above
equation , using 3.35 we get
\begin{eqnarray*}
F_{ji}^{\alpha }\left( b\right) &=&\left( id\otimes f_{ji}\right)
\Delta _{\alpha ,1}\left( b\right)
\\
&=&f_{ji}*b
\end{eqnarray*}
Inserting this result into 3.32 we obtain 3.28. In order to prove
3.29 we have to show that
$$\sum_{j}\left( f_{ji}*f_{hj}\circ
S_{1}^{-1}\right) =\delta _{ih}\varepsilon \eqno(3.37)$$ Let $a\in
A_{1}$ .Then
\begin{eqnarray*}
\sum_{j} \left( f_{ji}*f_{hj} \circ S_{1}^{-1} \right) \left(
S_{1}\left( a\right) \right) &=& \sum_{j}\left( f_{ji}\otimes
f_{hj} \circ S_{1}^{-1}\right) \Delta _{1,1}\left( S_{1}\left(
a\right) \right)
\\
&=&\sum_{j}\left( f_{ji}\otimes f_{hj} \circ S_{1}^{-1}\right)
\sigma _{A_{1},A_{1}}\left(S_{1}\otimes S_{1}\right) \Delta
_{1,1}\left( a\right)
\\
&=&\sum_{j}\left( f_{hj}\otimes f_{ji}\right) \left( id\otimes
S_{1}\right) \Delta_{1,1}\left( a\right)
\\
&=&\sum_{j}f_{hi}m_{1}\left(id\otimes S_{1}\right) \Delta
_{1,1}\left( a\right)
\\
&=&\sum_{j}f_{hi}\left( \varepsilon \left( a\right) 1_{1}\right)
\\
&=&\sum_{j}f_{hi}\left(\varepsilon \left( S_{1}\left( a\right)
\right) 1_{1}\right)
\\
&=&\sum_{j}f_{hi}\left(1_{1}\right) \varepsilon \left( S_{1}\left(
a\right) \right)
\\
&=&\delta _{hi}\varepsilon \left( S_{1}\left( a\right) \right)
\end{eqnarray*}
i.e.
$$\sum_{j}f_{ji}*\left( f_{hj}\circ
S_{1}^{-1}\right) =\delta _{ih}\varepsilon $$ Similarly , one can
check that
$$\sum_{j}\left( f_{jh}\circ
S_{1}^{-1}\right) *f_{ij}=\delta _{hi}\varepsilon \eqno(3.38)$$
From equation 3.28 we have that for any $\alpha \in \pi $ , $b\in
A_{\alpha } , \omega _{j}\in _{inv}\Gamma _{\alpha }$
$$\omega _{j}b=\sum_{h}\left(
f_{jh}*b\right) \omega _{h}$$ Inserting in this equation $b=\left(
f_{jh}\circ S_{1}^{-1}\right) *a$ for some $a\in A_{\alpha }$ and
summing over $j$ we obtain
\begin{eqnarray*}
\sum_{j}\omega _{j}\left( f_{jh}\circ S_{1}^{-1}\right) *a &=&
\sum_{j,h}\left( f_{jh}*\left( \left( f_{jh}\circ
S_{1}^{-1}\right) *a\right) \right) \omega _{h}
\\
&=&\sum_{j,h}\left( \left( f_{jh}*\left( f_{jh}\circ
S_{1}^{-1}\right) \right) *a\right) \omega _{h}
\\
&=&\sum_{j,h}\delta _{ih}\left( \varepsilon *a\right) \omega _{h}
\\
&=&a\omega _{i}
\end{eqnarray*}
Recall that $\varepsilon *a=\left( id\otimes \varepsilon \right)
\Delta _{\alpha ,1}\left( a\right) =a$ , and hence 3.29 follows.

To prove 2:
 For any $\alpha \in \pi ,\rho \in \Gamma _{\alpha
},$we have from statement 1 and formula 3.29 that
\begin{eqnarray*}
\rho &=&\sum_{i}a_{i}\omega _{i},\;\;\;\;\; a_{i}\in A_{\alpha
},\; \omega _{i}\in _{inv}\Gamma _{\alpha},\; i \in I
\\
&=&\sum_{i,j}\omega _{j}\left( \left( f_{ij}\circ
S_{1}^{-1}\right) *a_{i}\right)
\\
&=&\sum_{j}\omega _{j}b_{j},
\end{eqnarray*}
where

$$ b_{j}=\sum_{i}\left( f_{ij}\circ S_{1}^{-1}\right)
*a_{i}\in A_{\alpha },\;\;\forall j \in I.$$

For uniqueness:

Assume that for some $b_{i}$ $\left( i\in I\textrm{ only finite number of }%
b_{i}^{,}s\textrm{ are different from zero}\right) $we have:
$$\sum_{i}\omega_{i}b_{i}=0$$
We have to show that all $b_{i}'s$ $=0 \left( i\in I\textrm{
}\right) $ .Using the uniqueness of decomposition 3.26  we have
$$\sum_{i}\omega_{i}b_{i}=0 $$
Then
$$\sum_{i,j}\left(
f_{i,j}*b_{i}\right) \omega _{j}=0$$

$$\sum_{i}\left(
f_{i,j}*b_{i}\right) =0\qquad \;\;\; \forall j\in I $$

Computing the convolution product with $f_{jh}\circ S_{1}^{-1}$
summing over $j$ and using 3.38
\begin{eqnarray*}
0&=&\sum_{i,j}\left( f_{jh}\circ
S_{1}^{-1}\right) *\left( f_{i,j}*b_{i}\right) \\
&=&\sum_{i,j}\left( \left( f_{jh}\circ
S_{1}^{-1}\right) *f_{i,j}\right) *b_{i} \\
&=&\sum_{i,j}\delta _{hi}\left(
\varepsilon *b_{i}\right) \\
&=&b_{i}
\end{eqnarray*}
i.e. $b_{i}=0$ for each $i\in I.$
\end{proof}

Theorem 3.3 gives the complete description of left covariant $\pi
-$graded bimodules . Using 3.28 , 3.1 we have

$$\left( \sum_{i}a_{i}\omega _{i}\right)
b=\sum_{i}a_{i}\left( \omega _{i}b\right) =\sum_{i,.j}a_{i}\left(
f_{ij}*b\right) \omega _{j}\eqno(3.39)$$

$$\Delta _{\alpha ,\beta }^{l}\left(
\sum_{i}a_{i}\omega _{i}\right) =\sum_{i}\Delta _{\alpha ,\beta
}\left( a_{i}\right) \Delta _{\alpha ,\beta }^{l}\left( \omega
_{i}\right)=\sum_{i}\Delta _{\alpha ,\beta }\left( a_{i}\right)
\left( 1_{\alpha }\otimes \xi _{i}\right) ,\xi _{i}\in
_{inv}\Gamma _{\beta }\eqno(3.40)$$

If $\left( f_{ij}\right) _{i,j\in I}$ is a family of linear functionals in $%
A^{^{\prime }}=\oplus _{\alpha \in \pi }A_{\alpha }^{^{\prime }}$
satisfying relations 3.30, 3.31, then considering the left module
$\Gamma =\left\{ \Gamma _{\alpha }\right\}_{\alpha\in\pi} $
generated by $\omega _{i}^{\alpha },\alpha \in \pi, i\in I$ , and
using the above formulae to introduce the right multiplication by
elements of $A$ , and the left action of $A$ we obtain a left
covariant $\pi -$graded bimodule .

\begin{defn}
Let $\left( \Gamma ,\Delta ^{r}\right) $ be a right covariant $\pi
$-graded bimodule over A . An element $\eta \in \Gamma _{\alpha }$
is said to be right invariant if
$$\Delta
_{\alpha ,1}^{r}\left( \eta \right) =\eta \otimes
1_{1}\eqno(3.41)$$
\end{defn}
Denote by $\Gamma _{inv}=\left\{ \Gamma _{inv}^{\alpha }\right\} $
the set of all left invariant elements of $\Gamma $ .
Clearly,$\Gamma _{inv}^{\alpha }$ is a linear subspace of $\Gamma
_{\alpha }$ for each $\alpha \in \pi $.
\begin{thm}
Let $\Gamma =\left( \left\{ \Gamma ^{\alpha }\right\} _{\alpha \in
\pi
},\Delta ^{r}\right) $ be a right covariant $\pi -$graded bimodule over A,$%
\left\{ \eta _{i}^{\alpha }\right\} _{\alpha \in \pi }$be a basis
of $\Gamma _{inv}^{\alpha }$ of all right invariant elements of
$\Gamma _{\alpha }$ for each $\alpha \in \pi $ .Then
\begin{enumerate}
\item For any $\alpha \in \pi $ ,any element $\varrho \in $
$\Gamma _{\alpha }$ is of the form
$$\varrho
=\sum_{i}a_{i}\eta _{i}\eqno(3.42)$$

where $a_{i}$ $^{,}s\in A_{\alpha }$ are uniquely determined ,
$\eta _{i}^{,}s\in $ $\Gamma _{inv}^{\alpha }$,for any $\alpha \in
\pi $ . \item For any $\alpha \in \pi $ ,any element $\rho \in $
$\Gamma _{\alpha }$ is of the form
$$\varrho =\sum_{i}\eta
_{i}b_{i}\eqno(3.43)$$ where $b_{i}$ $^{,}s\in A_{\alpha }$ are
uniquely determined , $\eta _{i}^{,}s\in \Gamma _{inv}^{\alpha
}$,for any $\alpha \in \pi $ . \item There exists linear
functionals $g_{ij}$ $\in A^{^{\prime }}=\oplus _{\alpha \in \pi
}A_{\alpha }^{^{\prime }}$ such that for any $\alpha \in \pi $
$$\eta_{i}b=\sum_{j}\left( b*g_{ij}\right) \eta _{j}\eqno(3.44)$$
$$a\eta _{i}=\sum_{j}\eta
_{j}\left( a*\left( g_{ij}\circ S_{1}^{-1}\right) \right)
\eqno(3.45)$$ where $a,b\in A_{\alpha },\eta _{i}^{,}s,\eta
_{j}^{,}s\in \Gamma _{inv}^{\alpha }$ .These functionals are
uniquely determined by $3.44.$They satisfy the following relations
$$g_{ij}(ab)=\sum_{k}g_{ik}(a)g_{kj}(b)\eqno(3.46)$$
for any $i,j\in I$ , $a,b\in A_{\alpha }$ .Moreover
$$g_{ij}(1_{\alpha
})=\delta _{ij}\eqno(3.47)$$
\end{enumerate}
The proof is similar to that of theorem 3.3.
\end{thm}
\begin{rem}
Any functional $g_{ij}\in A^{^{\prime }}=\oplus _{\alpha \in \pi
}A_{\alpha }^{^{\prime }}$ is of the form $g_{ij}=\sum_{\alpha }$
$g_{ij}$ $^{\alpha }$ where
$$g_{ij}^{\alpha }\left(
a\right) =0 \;\;\;\;\; if \; a\notin A_{\alpha }$$
\end{rem}
%=======================================================================
\begin{thm}
Let $\Gamma =\left( \left\{ \Gamma _{\alpha }\right\} _{\alpha \in
\pi },\Delta ^{l},\Delta ^{r}\right) $ be a $\pi -$graded
bicovariant bimodule over A,$\left\{ \left( \omega _{i}^{\alpha
}\right)_{i\in \ I} \right\} _{\alpha \in \pi }$ be a basis of
$_{inv}\Gamma =\left\{ _{inv}\Gamma _{\alpha }\right\} _{\alpha
\in \pi }$ of all left invariant elements of $\Gamma $ .Then
\begin{enumerate}
\item For any $i\in I,\alpha ,\beta \in \pi ,\omega _{i}^{\alpha
\beta}\in \Gamma _{\alpha \beta }$
$$ \Delta _{\alpha ,\beta }^{r}\left( \omega_{i}^{\alpha \beta}\right) =\sum_{j}\omega _{j}^{\alpha }\otimes R_{ji} \eqno(3.48)$$
where $i,j\in\pi  , R_{ji}\in A_{\beta }$ satisfy the following
relation
$$\Delta _{\alpha ,\beta }\left( R_{ji}\right)
=\sum_{h}R_{jh}\otimes R_{hi} \eqno(3.49)$$ and for $R_{ji}\in
A_{1}$
$$ \varepsilon \left( R_{ji}\right) =\delta_{ji} \eqno(3.50) $$
\item For each $\alpha \in \pi $ there exists a basis $\left( \eta
_{i}\right) _{i\in I}$ of all right invariant elements of $\Gamma
_{\alpha }$ such that for $\omega_{i}\in\Gamma_{\alpha}$
$$\omega _{i}=\sum_{j}\eta _{j}R_{ji}\qquad \;\;\; \forall i\in I\eqno(3.51)$$
\item For any $j,h\in I$ $,a\in A_{\alpha }$
$$R_{ij}\left( a*f_{ih}\right)=\left(g_{ji}*a\right) R_{hi},\;\;\; i,j\in I \eqno(3.52)$$
where $f_{ij},g_{ij}$ are functionals introduced in theorems 3.3 ,
3.4
\end{enumerate}
\end{thm}
\begin{proof}
Using equation 3.9 for any $\alpha ,\beta ,\gamma \in \pi $ we
have
$$\left( \Delta _{\alpha ,\beta }^{l}\otimes
id\right) \Delta _{\alpha \beta ,\gamma }^{r}=\left( id\otimes
\Delta _{\beta ,\gamma }^{r}\right) \Delta _{\alpha ,\beta \gamma
}^{l}$$ Let $\omega _{i}^{\alpha \beta \gamma }\in \Gamma _{\alpha
\beta \gamma }$
\begin{eqnarray*}
\left( \Delta _{\alpha ,\beta }^{l}\otimes id\right) \Delta
_{\alpha \beta ,\gamma }^{r}\left( \omega _{i}^{\alpha \beta
\gamma }\right) &=&\left( id\otimes \Delta _{\beta ,\gamma
}^{r}\right) \Delta _{\alpha ,\beta \gamma }^{l}\left( \omega
_{i}^{\alpha \beta \gamma }\right)
\\
&=&1_{\alpha }\otimes \Delta _{\beta ,\gamma }^{r}\left( \omega
_{i}^{\beta \gamma }\right)
\end{eqnarray*}
i.e.
 $$\Delta _{\alpha \beta ,\gamma }^{r}\left(
\omega _{i}\right) \in _{inv}\Gamma _{\alpha \beta }\otimes
A_{\gamma }$$

Then for $\omega _{i}^{\alpha \beta \gamma }\in \Gamma _{\alpha
\beta \gamma }$

$$\Delta _{\alpha ,\beta \gamma }^{r}\left( \omega
_{i}^{\alpha \beta \gamma }\right) =\sum_{j}\omega _{j}^{\alpha
}\otimes R_{ji}$$

Applying $\left( id\otimes \Delta _{\beta ,\gamma }\right) $to
both sides of the above equation
\begin{eqnarray*}
\sum_{j}\omega _{j}^{\alpha }\otimes \Delta _{\beta ,\gamma
}\left( R_{ji}\right) &=&\left( id\otimes \Delta _{\beta ,\gamma
}\right) \Delta _{\alpha ,\beta \gamma }^{r}\left( \omega
_{i}^{\alpha \beta \gamma }\right)
\\
&=&\left( \Delta _{\alpha ,\beta }^{r}\otimes id\right) \Delta
_{\alpha \beta ,\gamma }^{r}\left( \omega _{i}^{\alpha \beta
\gamma }\right)
\\
&=&\left( \Delta _{\alpha ,\beta }^{r}\otimes id\right) \left(
\sum_{h}\omega _{h}^{\alpha \beta }\otimes R_{hi}\right)
\\
&=&\sum_{j,h}\omega _{j}^{\alpha }\otimes R_{jh}\otimes R_{hi}
\end{eqnarray*}
Comparing both sides of the above equation
$$\Delta _{\beta ,\gamma }\left( R_{ji}\right)
=\sum_{h}R_{jh}\otimes R_{hi}$$and hence 3.49 is proven.
 Let
$\omega _{i}^{\alpha }\in \Gamma _{\alpha }$
$$\Delta _{\alpha ,1}\left( \omega _{i}^{\alpha }\right)
=\sum_{j}\omega _{j}^{\alpha }\otimes R_{ji} \;\;  ,R_{ji}\in
A_{1}$$

Applying $\left( id\otimes \varepsilon \right) $ to both sides of
the above equation
\begin{eqnarray*}
\left( id\otimes \varepsilon \right) \Delta _{\alpha ,1}\left(
\omega _{i}^{\alpha }\right) &=&\omega _{i}^{\alpha }
\\
&=&\left( id\otimes \varepsilon \right) \left( \sum_{j}\omega
_{j}^{\alpha }\otimes R_{ji}\right)
\\
&=&\sum_{j}\omega _{j}^{\alpha }\otimes \varepsilon \left(
R_{ji}\right)
\end{eqnarray*}
$\Longrightarrow $
$$\varepsilon \left(
R_{ji}\right) =\delta _{ji}$$ To prove statement 2: First we have
that for $R_{ij}\in\ A_{1} , \alpha \in \pi $
%================================================================
$$m_{\alpha }\left( id\otimes S_{\alpha
^{-1}}\right) \Delta _{\alpha ,\alpha ^{-1}}=m_{\alpha }\left(
S_{\alpha ^{-1}}\otimes id\right) \Delta _{\alpha ^{-1},\alpha
}=\varepsilon 1_{\alpha }$$

By using 3.49 ,3.50  we obtain

$$\sum_{h}S_{\alpha ^{-1}}\left( R_{ih}\right)
R_{hj}=\delta _{ij}1_{\alpha }\eqno(3.53)$$

$$\sum_{h}R_{ih}S_{\alpha ^{-1}}\left(
R_{hj}\right) =\delta _{ij}1_{\alpha }\eqno(3.54)$$

For any $\alpha \in \pi ,j\in I,$ let

$$\eta _{j}=\sum_{i}\omega _{i}S_{\alpha
^{-1}}\left( R_{ij}\right) \eqno(3.55)$$

Multiplying both sides of  3.55  by $R_{ji}$ and summing over $j$
then using  3.53  we obtain
\begin{eqnarray*}
\sum_{j}\eta _{j}R_{ji}&=&\sum_{i,j}\omega _{i}S_{\alpha
^{-1}}\left( R_{ij}\right) R_{ji}
\\
&=&\omega _{i}
\end{eqnarray*}
and 3.51 follows. It remains to show that $\eta _{j}$ defined in
3.55 is right invariant

Let $\eta _{j}\in \Gamma _{\alpha },$ $\eta _{j}=\sum_{i}\omega
_{i}S_{\alpha ^{-1}}\left( R_{ij}\right) ,\omega _{i}\in
_{inv}\Gamma _{\alpha },R_{ij}\in A_{\alpha ^{-1}}.$
\begin{eqnarray*}
\Delta _{\alpha ,1}^{r}\left( \eta _{j}\right) &=&\Delta _{\alpha
,1}^{r}\left( \sum_{i}\omega _{i}S_{\alpha ^{-1}}\left(
R_{ij}\right) \right)
\\
&=&\sum_{i}\Delta _{\alpha ,1}^{r}\left( \omega _{i}\right) \Delta
_{\alpha ,1}\left( S_{\alpha ^{-1}}\left( R_{ij}\right) \right)
\\
&=&\sum_{i}\Delta _{\alpha ,1}^{r}\left( \omega _{i}\right) \sigma
_{A_{\alpha },A_{1}}\left( S_{1}\otimes S_{\alpha ^{-1}}\right)
\Delta _{1,\alpha ^{-1}}\left( R_{ij}\right)
\\
&=&\sum_{i,h,k}\left( \omega _{h}\otimes R_{hi}\right) \left(
S_{\alpha ^{-1}}\left( R_{kj}\right) \otimes S_{1}\left(
R_{ik}\right) \right)
\\
&=&\sum_{i,h,k}\omega _{h}S_{\alpha ^{-1}}\left( R_{kj}\right)
\otimes R_{hi}S_{1}\left( R_{ik}\right)
\\
&=&\sum_{h,k}\omega _{h}S_{\alpha ^{-1}}\left( R_{kj}\right)
\otimes \delta _{h,k}1_{1}
\\
&=&\sum_{k}\omega _{k}S_{\alpha ^{-1}}\left( R_{kj}\right) \otimes
1_{1}
\\
&=&\eta _{j}\otimes 1_{1}
\end{eqnarray*}
For any $\alpha \in \pi ,$let $\eta \in \Gamma _{\alpha }$ be a
right invariant element .According to theorem 3.3.2
\begin{eqnarray*}
\eta &=&\sum_{i}\omega _{i}c_{i}\;\;\;\;,c_{i}\in A_{\alpha }
\\
&=&\sum_{i,j}\eta _{j}R_{ji}c_{i}\;\;\;\;,R_{ji}\in A_{\alpha }
\end{eqnarray*}
then
$$ \eta =\sum_{j}\eta
_{j}b_{j},\;\;\;\;b_{j}\in A_{\alpha }\eqno(3.56)$$ If $\eta
=\sum_{i}\omega _{i}S_{\alpha ^{-1}}\left( R_{ij}\right),$then
using 3.55 we have
$$\sum_{i,j}\omega _{i}S_{\alpha
^{-1}}\left( R_{ij}\right) b_{j}=0$$ Using 3.3.1 we have
$$\sum_{j}S_{\alpha ^{-1}}\left(
R_{ij}\right) b_{j}=0 \;\;\;\; for \; each \; i\in I $$

Multiplying both sides of the above equation by $R_{ji}$ we obtain

$$b_{j}=0$$
for any $j\in I.$

This means that the decomposition 3.56 is unique.

Applying $\Delta _{1,\alpha }^{r}$ to both sides of 3.56 \qquad
\begin{eqnarray*}
\Delta _{1,\alpha }^{r}\left( \eta \right) &=&\Delta _{1,\alpha
}^{r}\left( \sum_{j}\eta _{j}b_{j}\right)
\\
\xi ^{1}\otimes 1_{\alpha }&=&\sum_{j}\left( \xi _{j}^{1}\otimes
1_{\alpha }\right) \Delta _{1,\alpha }\left( b_{j}\right)
\end{eqnarray*}
Comparing this formula with 3.56 we get that

$$\Delta _{1,\alpha }\left(
b_{j}\right) =b_{j\left( 1,1\right) }\otimes 1_{\alpha }$$

Applying $\varepsilon \otimes id$ we get that $%
b_{j}=\varepsilon \left( b_{j\left( 1,1\right) }\right) 1_{\alpha
}.$This way we proved that for any $\alpha \in \pi $ ,any $\eta
\in \Gamma _{inv}^{\alpha }$ is unique linear combination of $\eta
_{j}\left( j\in I\right) .$Therefore ,$\left( \eta _{j}\right)
_{j\in I}$ is a basis in $\Gamma _{\alpha }^{inv}$ and statement 2
is proven.

To prove statement 3 :

Using 3.45 we have for any $\alpha \in \pi , a\in A_{\alpha } ,
\eta _{j}\in \Gamma _{inv}^{\alpha }$

$$a\eta _{j}=\sum_{i}\eta _{i}\left( a*\left(
g_{ji}\circ S_{1}^{-1}\right) \right) $$

Using 3.55 we get

$$\sum_{i}a\omega _{i}S_{\alpha ^{-1}}\left( R_{ij}\right)
=\sum_{i,h}\omega _{h}S_{\alpha ^{-1}}\left( R_{hi}\right) \left(
a*\left( g_{ji}\circ S_{1}^{-1}\right) \right) $$

Using 3.31 we get

$$\sum_{i,h}\omega _{h}\left( \left( f_{ih}\circ
S_{1}^{-1}\right) *a\right) S_{\alpha ^{-1}}\left( R_{ij}\right)
=\sum_{i,h}\omega _{h}S_{\alpha ^{-1}}\left( R_{hi}\right) \left(
a*\left( g_{ji}\circ S_{1}^{-1}\right) \right) $$

Using 3.3.1 we get

$$ \sum_{i}\left( \left( f_{ih}\circ S_{1}^{-1}\right)
*a\right) S_{\alpha ^{-1}}\left( R_{ij}\right) =\sum_{i}S_{\alpha
^{-1}}\left( R_{hi}\right) \left( a*\left( g_{ji}\circ
S_{1}^{-1}\right) \right) $$

Applying $S_{\alpha }\left( S_{\alpha }=S_{\alpha
^{-1}}^{-1}\right) $ to both sides of this equation ,using that
$S_{\alpha }$ is antimultiplicative we get:

$$ \sum_{i}R_{ij} S_{\alpha }\left(\left( f_{ih}\circ
S_{1}^{-1}\right) *a\right) =\sum_{i}S_{\alpha }\left( a*\left(
g_{ji}\circ S_{1}^{-1}\right) \right) R_{hi}\eqno(3.57)$$

We compute
\begin{eqnarray*}
S_{\alpha }\left( \left( f_{ih}\circ S_{1}^{-1}\right) *a\right)
&=&S_{\alpha ^{-1}}^{-1}\left( id\otimes \left( f_{ih}\circ
S_{1}^{-1}\right) \right) \Delta _{\alpha ,1}\left( a\right)
\\
&=&\left( id\otimes f_{ih}\right) \left( S_{\alpha
^{-1}}^{-1}\otimes S_{1}^{-1}\right) \Delta _{\alpha ,1}\left(
a\right)
\\
&=&\left( f_{ih}\otimes id\right) \Delta _{1,\alpha ^{-1}}\left(
S_{\alpha ^{-1}}^{-1}\left( a\right) \right)
\\
&=&S_{\alpha ^{-1}}^{-1}\left( a\right) *f_{ih}
\\
S_{\alpha }\left( a*\left( g_{ji}\circ \ S_{1}^{-1}\right) \right)
&=&S_{\alpha ^{-1}}^{-1}\left( \left( g_{ji}\circ
S_{1}^{-1}\right) \otimes id\right) \Delta _{1,\alpha }\left(
a\right)
\\
&=&\left( g_{ji}\otimes id\right) \left( S_{1}^{-1}\otimes
S_{\alpha ^{-1}}^{-1}\right) \Delta _{1,\alpha }\left( a\right)
\\
&=&\left( id\otimes g_{ji}\right) \Delta _{\alpha ^{-1},1}\left(
S_{\alpha ^{-1}}^{-1}\left( a\right) \right)
\\
&=&g_{ji}*S_{\alpha ^{-1}}^{-1}\left( a\right)
\end{eqnarray*}
i.e.
$$ \sum_{i}R_{ij}\left( S_{\alpha
^{-1}}^{-1}\left( a\right) *f_{ih}\right) =\sum_{i}\left(
g_{ji}*S_{\alpha ^{-1}}^{-1}\left( a\right) \right) R_{hi}$$
Replacing $a$ by $S_{\alpha ^{-1}}^{-1}\left( a\right) $ we obtain
$$\sum_{i}R_{ij}\left( a*f_{ih}\right)
=\sum_{i}\left( g_{ji}*a\right) R_{hi}$$ And 3.52 follows. Note
that if $ a , R_{ij} , R_{hi} \in A_{1} , $ then applying
$\varepsilon $ to both sides of 3.52 and using 3.50 we obtain:
\begin{eqnarray*}
\varepsilon \left( \sum_{i}R_{ij}\left( a*f_{ih}\right) \right)
&=&\varepsilon \left( \sum_{i}\left( g_{ji}*a\right) R_{hi}\right)
\\
\sum_{i}\varepsilon \left( R_{ij}\right) \varepsilon \left(
a*f_{ih}\right) &=&\sum_{i}\varepsilon \left( g_{ji}*a\right)
\varepsilon \left( R_{hi}\right)
\\
\sum_{i}\delta _{ij}\varepsilon \left( a*f_{ih}\right)
&=&\sum_{i}\varepsilon \left( g_{ji}*a\right) \delta _{hi}
\end{eqnarray*}
But
\begin{eqnarray*}
\varepsilon \left( a*f_{ih}\right) &=&\left( f_{ih}\otimes
id\right) \Delta _{1,1}\left( a\right)
\\
&=&f_{ih}\left( a_{\left( 1,1\right) }\right) \varepsilon \left(
a_{\left( 2,1\right) }\right)
\\
&=&f_{ih}\left( a\right)
\end{eqnarray*}
Similarly: $\varepsilon \left( g_{ji}*a\right) =g_{ji}\left(
a\right) $ i.e. $f_{ij}\left( a\right) =g_{ij}\left( a\right)$,
for any $a\in A_{1}.$ From which we get that

$$\sum_{i}R_{ij}\left( a*f_{ih}\right)
=\sum_{i}\left( f_{ji}*a\right) R_{hi}\eqno(3.58)$$
\end{proof}.

For any $\alpha ,\beta \in \pi ,\eta _{j}\in \Gamma _{inv}^{\alpha
\beta },$ applying $\Delta _{\alpha ,\beta }^{l}$ to both sides of
equation 3.55 we obtain:
\begin{eqnarray*}
\Delta _{\alpha ,\beta }^{l}\left( \eta _{j}\right) &=&\Delta
_{\alpha ,\beta }^{l}\left( \sum_{h}\omega _{h}S_{\left( \alpha
\beta \right) ^{-1}}\left( R_{hj}\right) \right)
\\
&=&\sum_{h}\Delta _{\alpha ,\beta }^{l}\left( \omega _{h}\right)
\Delta _{\alpha ,\beta }\left( S_{\left( \alpha \beta \right)
^{-1}}\left( R_{hj}\right) \right)
\\
&=&\sum_{i,h}\left( 1_{\alpha }\otimes \omega _{h}^{\beta }\right)
\left( \sigma _{A_{\beta },A_{\alpha }}\left( S_{\beta
^{-1}}\otimes S_{\alpha ^{-1}}\right) \Delta _{\beta ^{-1},\alpha
^{-1}}\left( R_{hj}\right) \right)
\\
&=&\sum_{i,h}\left( 1_{\alpha }\otimes \omega _{h}^{\beta }\right)
\left( S_{\alpha ^{-1}}\left( R_{ij}\right) \otimes S_{\beta
^{-1}}\left( R_{hi}\right) \right)
\\
&=&\sum_{i,h}\left( S_{\alpha ^{-1}}\left( R_{ij}\right) \otimes
\omega _{h}^{\beta }S_{\beta ^{-1}}\left( R_{hi}\right) \right)
\end{eqnarray*}
i.e.
$$\Delta _{\alpha ,\beta }^{l}\left( \eta _{j}\right)
=\sum_{i}S_{\alpha ^{-1}}\left( R_{ij}\right) \otimes \eta
_{i}^{\beta }\eqno(3.59)$$ Using 3.5 ,3.48
$$\Delta _{\alpha ,\beta }^{r}\left(
\sum_{i}a_{i}\omega _{i}\right) =\sum_{i}\Delta _{\alpha ,\beta
}\left( a_{i}\right) \left( \omega _{j}^{\alpha}\otimes\
R_{ij}\right) \eqno(3.60)$$
%======================================================================
\begin{thm}
Let $\left( f_{ij}\right) _{i,j\in I}$ be the family of
functionals defined on $A$ satisfying relations 3.30, 3.31,
$\left( R_{ij}^{\alpha }\right) _{i,j\in I }$be a family of
elements of $A=\left\{ A_{\alpha }\right\} _{\alpha \in \pi }$
satisfying relations 3.49, 3.50, 3.58 for each $\alpha \in \pi .$
Consider the left module $\Gamma =\left\{ \Gamma _{\alpha
}\right\} _{\alpha \in \pi }$ over $A=\left\{ A_{\alpha }\right\}
_{\alpha \in \pi }$ generated by $\omega _{i}^{\alpha }$ $,$ $i\in
I,$ $\alpha \in \pi $ for each $\alpha \in \pi $ , and using
formulae 3.39, 3.40, 3.60 to introduce right multiplication by
elements of $A$, left and right actions of $A$ on $\Gamma $ then
$\Gamma =\left( \left\{ \Gamma _{\alpha }\right\} _{\alpha \in \pi
},\Delta ^{l},\Delta ^{r}\right) $ is a $\pi -$graded bicovariant
bimodule over $A.$
\end{thm}
\begin{proof} Using formula 3.39 to introduce right multiplication by
elements of $A$ ,one
can easily check that $\Gamma $ is also a $\pi -$graded right module over $%
A.$ , i.e.$\Gamma =\left\{ \Gamma _{\alpha }\right\} _{\alpha \in
\pi }$ is a $\pi -$graded bimodule over $A$ .

Using 3.40 to define a left action of $A$ on $\Gamma $ ,taking
into consideration 3.3.1 we find that 3.1 ,3.2 are satisfied for
let $\rho \in \Gamma _{\alpha \beta },$ $b\in A_{\alpha \beta
},\alpha ,\beta \in \pi ,$
using 3.3.1 $\rho =\sum_{i}a_{i}\omega _{i}^{\alpha \beta},$ $a_{i}\in A_{\alpha \beta },$ $%
\omega _{i}^{\alpha \beta}\in _{inv}\Gamma _{\alpha \beta }$
\begin{eqnarray*}
\Delta _{\alpha ,\beta }^{l}\left( b\rho \right) &=&\Delta
_{\alpha ,\beta }^{l}\left( b\sum_{i}a_{i}\omega _{i}^{\alpha
\beta}\right)
\\
&=&\Delta _{\alpha ,\beta }^{l}\left( \sum_{i}\left( ba_{i}\right)
\omega _{i}^{\alpha \beta}\right)
\\
&=&\sum_{i}\Delta _{\alpha ,\beta }\left( ba_{i}\right) \Delta
_{\alpha ,\beta }^{l}\left( \omega _{i}^{\alpha \beta}\right)
\\
&=&\sum_{i}\Delta _{\alpha ,\beta }\left( b\right) \Delta _{\alpha
,\beta }\left( a_{i}\right) \Delta _{\alpha ,\beta }^{l}\left(
\omega _{i}^{\alpha \beta}\right)
\\
&=&\Delta _{\alpha ,\beta }\left( b\right) \sum_{i}\Delta _{\alpha
,\beta }\left( a_{i}\right) \Delta _{\alpha ,\beta }^{l}\left(
\omega _{i}^{\alpha \beta}\right)
\\
&=&\Delta _{\alpha ,\beta }\left( b\right) \Delta _{\alpha ,\beta
}^{l}\left( \sum_{i}a_{i}\omega _{i}^{\alpha \beta}\right)
\\
&=&\Delta _{\alpha ,\beta }\left( b\right)\Delta _{\alpha ,\beta
}^{l}\left( \rho \right)
\end{eqnarray*}
And
\begin{eqnarray*}
\Delta _{\alpha ,\beta }^{l}\left( \rho b\right) &=&\Delta
_{\alpha ,\beta }^{l}\left( \left( \sum_{i}a_{i}\omega
_{i}^{\alpha \beta}\right) b\right)
\\
&=&\Delta _{\alpha ,\beta }^{l}\left( \sum_{i}a_{i}\left( \omega
_{i}^{\alpha \beta}b\right) \right)
\\
&=&\Delta _{\alpha ,\beta }^{l}\left( \sum_{i,j}a_{i}\left( \left(
f_{ij}*b\right) \omega _{j}^{\alpha \beta}\right) \right)
\\
&=&\sum_{i,j}\Delta _{\alpha ,\beta }^{l}\left( \left( a_{i}\left(
f_{ij}*b\right) \right) \omega _{j}^{\alpha \beta}\right)
\\
\end{eqnarray*}

\begin{eqnarray*}
&=&\sum_{i,j}\Delta _{\alpha ,\beta }\left( a_{i}\left(
f_{ij}*b\right) \right) \Delta _{\alpha ,\beta }^{l}\left( \omega
_{j}^{\alpha \beta}\right)
\\
&=&\sum_{i,j}\left( \Delta _{\alpha ,\beta }\left( a_{i}\right)
\Delta _{\alpha ,\beta }\left( f_{ij}*b\right) \right) \left(
1_{\alpha }\otimes \omega _{j}^{\beta }\right)
\\
&=&\sum_{i,j}\left( \Delta _{\alpha ,\beta }\left( a_{i}\right)
\Delta _{\alpha ,\beta }\left( b_{\left( 1,\alpha \beta \right)
}f_{ij}\left( b_{\left( 2,1\right) }\right) \right) \right) \left(
1_{\alpha }\otimes \omega _{j}^{\beta }\right)
\\
&=&\sum_{i,j}\left( \Delta _{\alpha ,\beta }\left( a_{i}\right)
\left( b_{\left( 1,\alpha \right) }\otimes f_{ij}\left( b_{\left(
3,1\right) }\right) b_{\left( 2,\beta \right) }\right) \right)
\left( 1_{\alpha }\otimes \omega _{j}^{\beta }\right)
\\
&=&\sum_{i,j}\Delta _{\alpha ,\beta }\left( a_{i}\right) \left(
b_{\left( 1,\alpha \right) }\otimes f_{ij}\left( b_{\left(
3,1\right) }\right) b_{\left( 2,\beta \right) }\omega _{j}^{\beta
}\right)
\\
&=&\sum_{i,j}\Delta _{\alpha ,\beta }\left( a_{i}\right) \left(
b_{\left( 1,\alpha \right) }\otimes \left( f_{ij}*b_{\left(
2,\beta \right) }\right) \omega _{j}^{\beta }\right)
\\
&=&\sum_{i}\Delta _{\alpha ,\beta }\left( a_{i}\right) \left(
b_{\left( 1,\alpha \right) }\otimes \omega _{i}^{\beta }b_{\left(
2,\beta \right) }\right)
\\
&=&\sum_{i}\Delta _{\alpha ,\beta }\left( a_{i}\right) \left(
1_{\alpha }\otimes \omega _{i}^{\beta }\right) \left( b_{\left(
1,\alpha \right) }\otimes b_{\left( 2,\beta \right) }\right)
\\
&=&\sum_{i}\Delta _{\alpha ,\beta }\left( a_{i}\right) \Delta
_{\alpha ,\beta }^{l}\left( \omega _{i}\right) \Delta _{\alpha
,\beta }\left( b\right)
\\
&=&\Delta _{\alpha ,\beta }^{l}\left( \left( \sum_{i}a_{i}\omega
_{i}\right) \right) \Delta _{\alpha ,\beta }\left( b\right)
\\
&=&\Delta _{\alpha ,\beta }^{l}\left( \rho \right) \Delta _{\alpha
,\beta }\left( b\right)
\end{eqnarray*}
Moreover , using 3.3.1 for any $\alpha $ , $\beta $ , $\gamma \in \pi $ $%
,\rho \in \Gamma _{\alpha \beta \gamma }$ , $\rho =\sum_{i}$
$a_{i}\omega _{i}$ where $a_{i}\in A_{\alpha \beta \gamma }$ ,
$\omega _{i}\in _{inv}\Gamma _{\alpha \beta \gamma }$ we have:
\begin{eqnarray*}
\left( \Delta _{\alpha ,\beta }\otimes id\right) \Delta _{\alpha
\beta ,\gamma }^{l}\left( \rho \right) &=&\left( \Delta _{\alpha
,\beta }\otimes id\right) \Delta _{\alpha \beta ,\gamma
}^{l}\left( \sum_{i}a_{i}\omega _{i}^{\alpha \beta \gamma}\right)
\\
&=&\sum_{i}\left( \Delta _{\alpha ,\beta }\otimes id\right) \Delta
_{\alpha \beta ,\gamma }\left( a_{i}\right) \Delta _{\alpha \beta
,\gamma }^{l}\left( \omega _{i}^{\alpha \beta \gamma}\right)
\\
&=&\sum_{i}\left( \Delta _{\alpha ,\beta }\otimes id\right) \left(
\Delta _{\alpha \beta ,\gamma }\left( a_{i}\right) \left(
1_{\alpha \beta }\otimes \omega _{i}^{\gamma }\right) \right)
\\
&=&\sum_{i}\left( \Delta _{\alpha ,\beta }\otimes id\right) \Delta
_{\alpha \beta ,\gamma }\left( a_{i}\right) \left( \Delta _{\alpha
,\beta }\otimes id\right) \left( 1_{\alpha \beta }\otimes \omega
_{i}^{\gamma }\right)
\\
&=&\sum_{i}\left( \left( \Delta _{\alpha ,\beta }\otimes id\right)
\Delta _{\alpha \beta ,\gamma }\left( a_{i}\right) \right) \left(
1_{\alpha }\otimes 1_{\beta }\otimes \omega _{i}^{\gamma }\right)\\
\left( id\otimes \Delta _{\beta ,\gamma }^{l}\right) \Delta
_{\alpha ,\beta \gamma }^{l}\left( \rho \right) &=&\left(
id\otimes \Delta _{\beta ,\gamma }^{l}\right) \Delta _{\alpha
,\beta \gamma }^{l}\left( \sum_{i}a_{i}\omega _{i}^{\alpha \beta
\gamma}\right)\\
\end{eqnarray*}

\begin{eqnarray*}
&=&\sum_{i}\left( id\otimes \Delta _{\beta ,\gamma }^{l}\right)
\left( \Delta _{\alpha ,\beta \gamma }\left( a_{i}\right) \Delta
_{\alpha ,\beta \gamma }^{l}\left( \omega _{i}^{\alpha \beta
\gamma}\right) \right)\\
&=&\sum_{i}\left( id\otimes \Delta _{\beta ,\gamma }\right) \left(
\Delta _{\alpha ,\beta \gamma }\left( a_{i}\right) \right) \left(
id\otimes \Delta _{\beta ,\gamma }^{l}\right) \left( \Delta
_{\alpha ,\beta \gamma }^{l}\left( \omega _{i}^{\alpha \beta
\gamma}\right) \right)\\
&=&\sum_{i}\left( id\otimes \Delta _{\beta ,\gamma }\right) \left(
\Delta _{\alpha ,\beta \gamma }\left( a_{i}\right) \right) \left(
id\otimes \Delta _{\beta ,\gamma }^{l}\right) \left( 1_{\alpha
}\otimes \omega _{i}^{\beta \gamma }\right)\\
&=&\sum_{i}\left( id\otimes \Delta _{\beta ,\gamma }\right) \left(
\Delta _{\alpha ,\beta \gamma }\left( a_{i}\right) \right) \left(
1_{\alpha }\otimes 1_{\beta }\otimes \omega _{i}^{\gamma
}\right)\\
&=&\sum_{i}\left( \Delta _{\alpha ,\beta }\otimes id\right) \Delta
_{\alpha \beta ,\gamma }\left( a_{i}\right) \left( 1_{\alpha
}\otimes 1_{\beta }\otimes \omega _{i}^{\gamma }\right)
\end{eqnarray*}
i.e.
$$ \left( \Delta _{\alpha ,\beta }\otimes
id\right) \Delta _{\alpha \beta ,\gamma }^{l}= \left( id\otimes
\Delta _{\beta ,\gamma }^{l}\right) \Delta _{\alpha ,\beta \gamma
}^{l}$$
%=======================================================================
which means that 3.3 is satisfied.

Finally , for any  $\alpha \in \pi ,$ letting  $\rho \in \Gamma
_{\alpha }$ , then using 3.3.1 $\rho =\sum_{i}a_{i}\omega
_{i}^{\alpha },a_{i}\in A_{\alpha }, \omega _{i}\in _{inv}\Gamma
_{\alpha }$
\begin{eqnarray*} \left(
\varepsilon \otimes id\right) \Delta _{1,\alpha }^{l}\left( \rho
\right) &=&\left( \varepsilon \otimes id\right) \Delta _{1,\alpha
}^{l}\left(
\sum_{i}a_{i}\omega _{i}\right)\\
\end{eqnarray*}

\begin{eqnarray*}
&=&\sum_{i}\left( \varepsilon \otimes id\right) \Delta _{1,\alpha
}\left( a_{i}\right) \Delta _{1,\alpha
}^{l}\left( \omega _{i}\right)\\
&=&\sum_{i}\left( \varepsilon \otimes id\right) \Delta _{1,\alpha
}\left( a_{i}\right) \Delta _{1,\alpha
}^{l}\left( \omega _{i}\right)\\
&=&\sum_{i}\left( \varepsilon \otimes id\right) \left( \Delta
_{1,\alpha }\left( a_{i}\right) \right) \left( \varepsilon \otimes
id\right) \left( \Delta
_{1,\alpha }^{l}\left( \omega _{i}\right) \right)\\
&=&\sum_{i}a_{i}\omega _{i}\\
&=&\rho
\end{eqnarray*}
i.e. 3.4 is satisfied.\\
 This means that $\Gamma =\left( \left\{
\Gamma _{\alpha }\right\} _{\alpha \in \pi },\Delta ^{l}\right) $
is a $\pi -$graded left
covariant bimodule over $A.$\\
Using formula 3.60 to introduce right action of $A$ on $\Gamma $
one can easily check that $\Gamma =\left( \left\{ \Gamma _{\alpha
}\right\} _{\alpha \in \pi },\Delta ^{r}\right) $ is a $\pi
-$graded right covariant bimodule over $A,$ for let $\rho \in
\Gamma _{\alpha \beta },$ $b\in A_{\alpha \beta },\alpha ,\beta
\in \pi ,$ using  3.3.1 $\rho =\sum_{i}a_{i}\omega _{i},$
$a_{i}\in A_{\alpha \beta },$ $\omega _{i}\in _{inv}\Gamma
_{\alpha \beta }$
\begin{eqnarray*}
\Delta _{\alpha ,\beta }^{r}\left( b\rho \right) &=& \Delta
_{\alpha ,\beta }^{r}\left( b\sum_{i}a_{i}\omega _{i}\right)\\
&=&\Delta _{\alpha ,\beta }^{r}\left( \sum_{i}\left( ba_{i}\right)\omega _{i}\right)\\
&=&\sum_{i} \Delta _{\alpha ,\beta }\left(ba_{i}\right) \Delta _{\alpha ,\beta }^{r}\left( \omega_{i}\right) \\
&=&\sum_{i}\Delta _{\alpha ,\beta }\left( b\right) \Delta _{\alpha
,\beta }\left( a_{i}\right) \Delta _{\alpha ,\beta}^{r}\left( \omega _{i}\right) \\
&=&\Delta _{\alpha ,\beta }\left( b\right) \sum_{i}\Delta
_{\alpha,\beta }\left( a_{i}\right) \Delta _{\alpha ,\beta}^{r}\left(\omega _{i}\right) \\
&=&\Delta _{\alpha ,\beta }\left( b\right) \Delta _{\alpha ,\beta
}^{r}\left( \sum_{i}a_{i}\omega _{i}\right) \\
&=&\Delta _{\alpha,\beta }\left( b\right) \Delta _{\alpha
,\beta }^{r}\left( \rho\right)\\
\end{eqnarray*}
Again , for $\rho \in \Gamma _{\alpha \beta },$ $b\in A_{\alpha
\beta },\alpha ,\beta \in \pi ,$ using  3.3.1 $\rho
=\sum_{i}a_{i}\omega _{i},$ $a_{i}\in A_{\alpha \beta },$ $\omega
_{i}\in _{inv}\Gamma _{\alpha \beta }.$ Using  3.44  we get
\begin{eqnarray*}
\Delta _{\alpha ,\beta }^{r}\left( \rho b\right) &=&\Delta
_{\alpha ,\beta }^{r}\left( \left( \sum_{i}a_{i}\omega _{i}\right)
b\right) \\
&=&\Delta _{\alpha ,\beta }^{r}\left( \sum_{i}a_{i}\left( \omega
_{i}b\right) \right) \\
&=&\Delta _{\alpha ,\beta }^{r}\left( \sum_{i,j}a_{i}\left( \left(
f_{ij}*b\right) \omega _{j}\right) \right) \\
&=&\Delta _{\alpha ,\beta }^{r}\left( \sum_{i,j}\left( a_{i}\left(
f_{ij}*b\right) \right) \omega _{j}\right) \\
&=&\sum_{i,j}\Delta _{\alpha ,\beta }\left( a_{i}\left(
f_{ij}*b\right) \right) \Delta _{\alpha ,\beta }^{r}\left( \omega
_{j}\right) \\
&=&\sum_{i,j}\Delta _{\alpha ,\beta }\left( a_{i}\right) \Delta
_{\alpha ,\beta }\left( f_{ij}*b\right) \Delta _{\alpha ,\beta
}^{r}\left( \omega _{j}\right) \\
&=&\sum_{i,j}\Delta _{\alpha ,\beta }\left( a_{i}\right) \Delta
_{\alpha ,\beta }\left( b_{\left( 1,\alpha \beta \right)
}f_{ij}\left( b_{2,1}\right) \right) \left( \omega _{k}\otimes
R_{kj}\right) \\
&=&\sum_{i,j,k}\Delta _{\alpha ,\beta }\left( a_{i}\right) \left(
b_{\left( 1,\alpha \right) }\otimes b_{\left( 2,\beta \right)
}f_{ij}\left( b_{3,1}\right) \right) \left( \omega _{k}\otimes
R_{kj}\right)\\
&=&\sum_{i,j,k}\Delta _{\alpha ,\beta }\left( a_{i}\right) \left(
b_{\left( 1,\alpha \right) }\omega _{k}\otimes b_{\left( 2,\beta
\right) }f_{ij}\left( b_{3,1}\right) R_{kj}\right) \\
\end{eqnarray*}

\begin{eqnarray*}
&=&\sum_{i,j,k}\Delta _{\alpha ,\beta }\left( a_{i}\right) \left(
b_{\left( 1,\alpha \right) }\omega _{k}\otimes \left(
f_{ij}*b_{\left( 2,\beta \right) }\right) R_{kj}\right) \\
&=&\sum_{i,j,k,l}\Delta _{\alpha ,\beta }\left( a_{i}\right)
\left( \omega _{l}\left( \left( f_{kl}\circ S_{1}^{-1}\right)
*b_{\left( 1,\alpha \right) }\right) \otimes R_{ji}\left(
b_{\left( 2,\beta \right) }*f_{jk}\right) \right) \\
&=&\sum_{i,j,k,l}\Delta _{\alpha ,\beta }\left( a_{i}\right)
\left( \left( \omega _{l}\otimes R_{ji}\right) \left( \left(
f_{kl}\circ S_{1}^{-1}\right) *b_{\left( 1,\alpha \right) }\otimes
b_{\left( 2,\beta \right) }*f_{jk}\right) \right) \\
&=&\sum_{i,j,k,l}\Delta _{\alpha ,\beta }\left( a_{i}\right)
\left( \left( \omega _{l}\otimes R_{ji}\right) \left( b_{\left(
1,\alpha \right) }f_{kl}\left( S_{1}^{-1}\left( b_{\left(
2,1\right) }\right) \right) \otimes f_{jk}\left( b_{\left(
3,1\right) }\right) b_{\left( 4,\beta \right) }\right) \right) \\
&=&\sum_{i,j,l}\Delta _{\alpha ,\beta }\left( a_{i}\right) \left(
\left( \omega _{l}\otimes R_{ji}\right) \left( b_{\left( 1,\alpha
\right) }\otimes f_{jl}\left( b_{\left( 3,1\right)
}S_{1}^{-1}\left( b_{\left( 2,1\right) }\right) \right) b_{\left(
4,\beta \right) }\right) \right) \\
&=&\sum_{i,j,l}\Delta _{\alpha ,\beta }\left( a_{i}\right) \left(
\left( \omega _{l}\otimes R_{ji}\right) \left( b_{\left( 1,\alpha
\right) }\otimes f_{jl}\left( \varepsilon \left( b_{\left(
2,1\right) }\right) 1_{1}\right) b_{\left( 3,\beta \right)
}\right) \right) \\
&=&\sum_{i,j,l}\Delta _{\alpha ,\beta }\left( a_{i}\right) \left(
\left( \omega _{l}\otimes R_{ji}\right) \left( b_{\left( 1,\alpha
\right) }\otimes \varepsilon \left( b_{\left( 2,1\right) }\right)
\delta _{jl}b_{\left( 3,\beta \right) }\right) \right) \\
&=&\sum_{i,j}\Delta _{\alpha ,\beta }\left( a_{i}\right) \left(
\left( \omega _{j}\otimes R_{ji}\right) \left( b_{\left( 1,\alpha
\right) }\otimes b_{\left( 2,\beta \right) }\right) \right)\\
&=&\sum_{i}\Delta _{\alpha ,\beta }\left( a_{i}\right) \Delta
_{\alpha ,\beta }^{r}\left( \omega _{i}\right) \Delta _{\alpha
,\beta }\left( b\right)\\
&=&\Delta _{\alpha ,\beta }^{r}\left( \rho \right) \Delta _{\alpha
,\beta }\left( b\right)\\
\end{eqnarray*}
and thus 3.5, 3.6 are satisfied.
Moreover , using 3.3.1 for any $\alpha $ , $\beta $ , $\gamma \in \pi $ $%
,\rho \in \Gamma _{\alpha \beta \gamma }$ , $\rho =\sum_{i}
a_{i}\omega _{i}^{\alpha \beta \gamma}$ where $a_{i}\in A_{\alpha
\beta \gamma }$ , $\omega _{i}^{\alpha \beta \gamma}\in
_{inv}\Gamma _{\alpha \beta \gamma }$ we have:
\begin{eqnarray*}
\left( \Delta _{\alpha ,\beta }^{r}\otimes id\right) \Delta
_{\alpha \beta ,\gamma }^{r}\left( \rho \right)&=&\left( \Delta
_{\alpha ,\beta }^{r}\otimes id\right) \Delta _{\alpha \beta
,\gamma }^{r}\left( \sum_{i}a_{i}\omega
_{i}\right)\\
&=&\sum_{i}\left( \Delta _{\alpha ,\beta }^{r}\otimes id\right)
\left( \Delta _{\alpha \beta ,\gamma }\left( a_{i}\right) \Delta
_{\alpha \beta ,\gamma }^{r}\left(
\omega _{i}\right) \right) \\
&=&\sum_{i}\left( \Delta _{\alpha ,\beta }\otimes id\right) \left(
\Delta _{\alpha \beta ,\gamma }\left( a_{i}\right) \right) \left(
\Delta _{\alpha ,\beta }^{r}\otimes id\right) \left( \Delta
_{\alpha \beta ,\gamma
}^{r}\left( \omega _{i}\right) \right)\\
&=&\sum_{i,j}\left( \Delta _{\alpha ,\beta }\otimes id\right)
\left( \Delta _{\alpha \beta ,\gamma }\left( a_{i}\right) \right)
\left( \Delta _{\alpha ,\beta }^{r}\otimes id\right) \left( \omega
_{j}^{\alpha \beta }\otimes R_{ji}\right)\\
&=&\sum_{i,j,k}\left( \Delta _{\alpha ,\beta }\otimes id\right)
\left( \Delta _{\alpha \beta ,\gamma }\left( a_{i}\right) \right)
\left( \omega _{k}^{\alpha }\otimes R_{kj}\otimes R_{ji}\right)\\
\left( id\otimes \Delta _{\beta ,\gamma }\right) \Delta _{\alpha
,\beta \gamma }^{r}\left( \rho \right) &=&\left( id\otimes\Delta
_{\beta ,\gamma }\right) \Delta _{\alpha ,\beta \gamma }^{r}\left(
\sum_{i}a_{i}\omega _{i}\right) \\
&=&\sum_{i}\left( id\otimes \Delta _{\beta ,\gamma }\right) \left(
\Delta _{\alpha ,\beta \gamma }\left( a_{i}\right) \Delta _{\alpha
,\beta \gamma }^{r}\left( \omega _{i}\right)
\right) \\
\end{eqnarray*}

\begin{eqnarray*}
&=&\sum_{i}\left( id\otimes \Delta _{\beta ,\gamma }\right) \left(
\Delta _{\alpha ,\beta \gamma }\left( a_{i}\right) \right) \left(
id\otimes \Delta _{\beta ,\gamma }\right) \left( \Delta _{\alpha
,\beta \gamma }^{r}\left(
\omega _{i}\right) \right) \\
&=&\sum_{i,k}\left( id\otimes \Delta _{\beta ,\gamma }\right)
\left( \Delta _{\alpha ,\beta \gamma }\left( a_{i}\right)
\right) \\
&=&\sum_{i,j,k}\left( id\otimes \Delta _{\beta ,\gamma }\right)
\left( \Delta _{\alpha ,\beta \gamma
}\left(a_{i}\right) \right) \left( \omega _{k}^{\alpha }\otimes R_{kj}\otimes R_{ji}\right) \\
\end{eqnarray*}
i.e. $\left( \Delta _{\alpha ,\beta }^{r}\otimes id\right) \Delta
_{\alpha \beta ,\gamma }^{r}=\left( id\otimes \Delta _{\beta
,\gamma }\right) \Delta _{\alpha ,\beta \gamma }^{r}$ which means
that 3.7 is satisfied.
 Finally, for any $\alpha \in \pi,$ letting $\rho \in \Gamma _{\alpha }$, then using
3.4.1  $\rho \in \Gamma _{\alpha }, \rho =\sum_{i} a_{i}\omega
_{i}$ where $a_{i}\in A_{\alpha }$ , $\omega _{i}\in _{inv}\Gamma
_{\alpha }$
\begin{eqnarray*}
\left( id\otimes \varepsilon \right) \Delta _{\alpha ,1}^{r}\left(
\rho \right) &=&\left( id\otimes \varepsilon \right) \Delta
_{\alpha ,1}^{r}\left(
\sum_{i}a_{i}\omega _{i}\right) \\
&=&\sum_{i}\left( id\otimes \varepsilon \right) \left( \Delta
_{\alpha ,1}\left( a_{i}\right) \Delta
_{\alpha ,1}^{r}\left( \omega _{i}\right) \right) \\
&=&\sum_{i}\left( id\otimes \varepsilon \right) \left( \Delta
_{\alpha ,1}\left( a_{i}\right) \right) \left( id\otimes
\varepsilon \right) \left( \Delta
_{\alpha ,1}^{r}\left( \omega _{i}\right) \right) \\
&=&\sum_{i}a_{i}\omega _{i}\\
&=&\rho \\
\end{eqnarray*}
 i.e. 3.8 is satisfied . This means that $\Gamma=\left( \left\{
\Gamma_{\alpha }\right\} _{\alpha \in \pi },\Delta ^{r}\right) $ is a $\pi -$%
graded right covariant bimodule over $A.$ To prove the
bicovariance conditions , for any $\alpha ,\beta ,\gamma \in \pi $
,$\rho \in \Gamma _{\alpha \beta \gamma }$, using 3.3.1  $\rho
=\sum_{i}a_{i}\omega _{i}^{\alpha \beta \gamma},$ $a_{i}\in
A_{\alpha \beta \gamma },$ $\omega _{i}^{\alpha \beta \gamma}\in
_{inv}\Gamma _{\alpha \beta \gamma }$  we compute
\begin{eqnarray*}
\left( id\otimes \Delta _{\beta ,\gamma }^{r}\right) \Delta
_{\alpha ,\beta \gamma }^{l}\left( \rho \right) &=&\left(
id\otimes \Delta _{\beta ,\gamma }^{r}\right) \Delta _{\alpha
,\beta \gamma }^{l}\left( \sum_{i}a_{i}\omega
_{i}^{\alpha \beta \gamma}\right)\\
&=&\sum_{i}\left( id\otimes \Delta _{\beta ,\gamma }^{r}\right)
\left( \Delta _{\alpha ,\beta \gamma }\left( a_{i}\right) \Delta
_{\alpha ,\beta \gamma }^{l}\left( \omega
_{i}^{\alpha \beta \gamma}\right) \right)\\
&=&\sum_{i}\left(id\otimes\Delta _{\beta ,\gamma }\right) \left(
\Delta _{\alpha ,\beta \gamma }\left( a_{i}\right) \right) \left(
id\otimes \Delta _{\beta ,\gamma }^{r}\right) \left( \Delta
_{\alpha ,\beta \gamma
}^{l}\left( \omega _.{i}^{\alpha \beta \gamma}\right) \right)\\
&=&\sum_{i}\left(id\otimes \Delta _{\beta ,\gamma }\right) \left(
\Delta _{\alpha ,\beta \gamma }\left( a_{i}\right) \right) \left(
1_{\alpha }\otimes \Delta _{\beta
,\gamma }^{r}\left( \omega _{i}^{\beta \gamma }\right) \right)\\
&=&\sum_{i,j}\left(id\otimes \Delta _{\beta ,\gamma }\right)
\left( \Delta _{\alpha ,\beta \gamma }\left( a_{i}\right) \right)
\left( 1_{\alpha }\otimes \omega _{j}^{\beta }\otimes
R_{ji}\right)\\
\end{eqnarray*}
\begin{eqnarray*}
  \left( \Delta _{\alpha ,\beta
}^{l}\otimes id\right) \Delta _{\alpha \beta ,\gamma }^{r}\left(
\rho \right)&=&\left( \Delta _{\alpha ,\beta }^{l}\otimes
id\right) \Delta _{\alpha \beta ,\gamma
}^{r}\left( \sum_{i}a_{i}\omega _{i}^{\alpha \beta \gamma}\right) \\
&=&\sum_{i}\left( \Delta _{\alpha ,\beta }^{l}\otimes id\right)
\left( \Delta _{\alpha \beta ,\gamma }\left( a_{i}\right) \Delta
_{\alpha \beta ,\gamma }^{r}\left(
\omega _{i}^{\alpha \beta \gamma}\right) \right) \\
&=&\sum_{i}\left( \Delta _{\alpha ,\beta }\otimes id\right) \left(
\Delta _{\alpha \beta ,\gamma }\left( a_{i}\right) \right) \left(
\Delta _{\alpha ,\beta }^{l}\otimes id\right) \left( \Delta
_{\alpha \beta ,\gamma
}^{r}\left( \omega _{i}^{\alpha \beta \gamma}\right) \right) \\
&=&\sum_{i,j}\left( \Delta _{\alpha ,\beta }\otimes id\right)
\left( \Delta _{\alpha \beta ,\gamma }\left( a_{i}\right) \right)
\left( \Delta _{\alpha ,\beta }^{l}\otimes id\right) \left( \omega
_{j}^{\alpha \beta }\otimes
R_{ji}\right)\\
&=&\sum_{i,j}\left( \Delta _{\alpha ,\beta }\otimes id\right)
\left( \Delta _{\alpha \beta ,\gamma }\left( a_{i}\right) \right)
\left( 1_{\alpha }\otimes \omega _{j}^{\beta
}\otimes R_{ji}\right) \\
\end{eqnarray*}
and hence 3.9 is proved and $\Gamma =\left( \left\{ \Gamma
_{\alpha
}\right\} _{\alpha \in \pi },\Delta ^{l},\Delta ^{r}\right) $is a $\pi -$%
graded bicovariant bimodule over $A.$
\end{proof}

\section{first order differential calculus on Hopf Group Coalgebras}
Let $A^{2}=\left\{ A_{\alpha }^{2}\right\} _{\alpha \in \pi }$ be the $\pi -$%
graded bimodule introduced in section 2.We introduce left and
right actions of $A$ on $A^{2}$ .For any $\alpha ,\beta \in \pi $
let $q\in A_{\alpha \beta }\otimes A_{\alpha \beta },$and $\left(
\Delta _{\alpha ,\beta }\otimes \Delta _{\alpha ,\beta }\right)
\left( q\right) =\sum_{k}a_{k}\otimes b_{k}\otimes c_{k}\otimes
d_{k}$, where $
a_{k},c_{k}\in A_{\alpha },$ $b_{k},d_{k}$ $\in A_{\beta },$ $k=1,2,...,n.$%
We set
$$\Phi _{\alpha ,\beta }^{l}\left(
q\right)=\sum_{k}a_{k}c_{k}\otimes b_{k}\otimes d_{k}\eqno(4.1)$$
$$\Phi _{\alpha ,\beta }^{r}\left( q\right) =\sum_{k}a_{k}\otimes
c_{k}\otimes b_{k}d_{k}\eqno(4.2)$$ We compute
\begin{eqnarray*}
\left( id\otimes m_{\beta }\right) \left( \Phi _{\alpha ,\beta
}^{l}\left( q\right) \right) &=&\left( id\otimes m_{\beta }\right)
\left(
\sum_{k}a_{k}c_{k}\otimes b_{k}\otimes d_{k}\right) \\
&=&\sum_{k}a_{k}c_{k}\otimes b_{k}d_{k}\\
 &=&\Delta _{\alpha,\beta }\left( m_{\alpha \beta }\left( q\right) \right) \\
 &=&0\\
\end{eqnarray*}
Similarly we have
\begin{eqnarray*}
\left( m_{\alpha }\otimes id\right) \left( \Phi _{\alpha ,\beta
}^{r}\left( q\right) \right) &=&\left( m_{\alpha }\otimes
id\right) \left( \sum_{k}a_{k}\otimes
c_{k}\otimes b_{k}d_{k}\right) \\
 \end{eqnarray*}

\begin{eqnarray*}
 &=&\sum_{k}a_{k}c_{k}\otimes b_{k}d_{k}\\
&=&\Delta _{\alpha ,\beta }\left( m_{\alpha \beta }\left( q\right)
\right) \\
& =0&
\end{eqnarray*}
Therfore,
$$\Phi _{\alpha ,\beta }^{l}:A_{\alpha \beta
}^{2}\longrightarrow A_{\alpha }\otimes A_{\beta }^{2}\eqno(4.3)$$
and
$$\Phi _{\alpha ,\beta }^{r}:A_{\alpha \beta
}^{2}\longrightarrow A_{\alpha }^{2}\otimes A_{\beta }\eqno(4.4)$$
Clearly, both are linear map. We will show that $A^{2}=\left(
\left\{ A_{\alpha }^{2}\right\} _{\alpha \in \pi },\Phi ^{l},\Phi
^{r}\right) $ is a $\pi -$ graded bicovariant bimodule over $A.$

First , we will prove that $A^{2}=\left( \left\{ A_{\alpha
}^{2}\right\} _{\alpha \in \pi },\Phi ^{l},\Phi ^{r}\right) $ is a
$\pi -$ graded left covariant bimodule over $A.$

Let $\alpha ,\beta \in \pi ,q\in A_{\alpha \beta }^{2} ,
q=b\otimes c,$ then
\begin{eqnarray*}
\Phi _{\alpha ,\beta }^{l}\left( aq\right) &=&\Phi _{\alpha
,\beta }^{l}\left( ab\otimes c\right) \\
&=&a_{\left( 1,\alpha \right) }b_{\left( 1,\alpha \right)
}c_{\left( 1,\alpha \right) }\otimes a_{\left( 2,\beta \right)
}b_{\left( 2,\beta \right) }\otimes c_{\left( 2,\beta \right) }\\
&=&\left( a_{\left( 1,\alpha \right) }\otimes a_{\left( 2,\beta
\right) }\right) \cdot \left( b_{\left( 1,\alpha \right)
}c_{\left( 1,\alpha \right) }\otimes b_{\left( 2,\beta \right)
}\otimes c_{\left( 2,\beta \right) }\right) \\
&=&\Delta _{\alpha ,\beta }\left( a\right) \cdot \Phi _{\alpha
,\beta }^{l}\left( b\otimes c\right) \\
&=&\Delta _{\alpha ,\beta }\left( a\right) \cdot \Phi _{\alpha
,\beta }^{l}\left( q\right) \\
\end{eqnarray*}
and
\begin{eqnarray*}
\Phi _{\alpha ,\beta }^{l}\left( qa\right) &=&\Phi _{\alpha ,\beta
}^{l}\left( b\otimes ca\right) \\
&=&b_{\left( 1,\alpha \right) }c_{\left( 1,\alpha \right)
}a_{\left( 1,\alpha \right) }\otimes b_{\left( 2,\beta \right)
}\otimes c_{\left( 2,\beta \right) }a_{\left( 2,\beta \right) }\\
&=&\left( b_{\left( 1,\alpha \right) }c_{\left( 1,\alpha \right)
}\otimes b_{\left( 2,\beta \right) }\otimes c_{\left( 2,\beta
\right) }\right) .\left( a_{\left( 1,\alpha \right) }\otimes
a_{\left( 2,\beta \right) }\right) \\
&=&\Phi _{\alpha ,\beta }^{l}\left( b\otimes c\right) \cdot \Delta
_{\alpha ,\beta }\left( a\right) \\
&=&\Phi _{\alpha ,\beta }^{l}\left( q\right) \cdot \Delta _{\alpha
,\beta }\left( a\right) \\
\end{eqnarray*}
Moreover , for any $\alpha ,\beta ,\gamma \in \pi , q\in A_{\alpha
\beta \gamma }^{2},q=a\otimes b$  we compute
\begin{eqnarray*}
\left( \Delta _{\alpha ,\beta }\otimes id\right) \Phi _{\alpha
\beta ,\gamma }^{l}\left( q\right)& =&\left( \Delta _{\alpha
,\beta }\otimes id\right) \Phi _{\alpha
\beta ,\gamma }^{l}\left( a\otimes b\right) \\
&=&\left( \Delta _{\alpha ,\beta }\otimes id\right) \left(
a_{\left( 1,\alpha \beta \right) }b_{\left( 1,\alpha \beta \right)
}\otimes a_{\left( 2,\gamma \right)
}\otimes b_{\left( 2,\gamma \right) }\right) \\
&=&a_{\left( 1,\alpha \right) }b_{\left( 1,\alpha \right) }\otimes
a_{\left( 2,\beta \right) }b_{\left( 2,\beta \right) }\otimes
a_{\left( 3,\gamma \right) }\otimes b_{\left( 3,\gamma \right) }\\
\end{eqnarray*}
\begin{eqnarray*}
\left( id\otimes \Phi _{\beta ,\gamma }^{l}\right) \Phi _{\alpha
,\beta \gamma }^{l}\left( q\right)& =&\left( id\otimes \Phi
_{\beta ,\gamma }^{l}\right) \Phi
_{\alpha ,\beta \gamma }^{l}\left( a\otimes b\right)\\
&=&\left(id\otimes \Phi _{\beta ,\gamma }^{l}\right) \left(
a_{\left( 1,\alpha \right) }b_{\left( 1,\alpha \right) }\otimes
a_{\left( 2,\beta \gamma \right) }\otimes b_{\left(
2,\beta \gamma \right) }\right) \\
&=&a_{\left( 1,\alpha \right) }b_{\left( 1,\alpha \right) }\otimes
a_{\left( 2,\beta \right) }b_{\left( 2,\beta \right) }\otimes
a_{\left( 3,\gamma \right) }\otimes b_{\left( 3,\gamma \right) }\\
\end{eqnarray*}
i.e. $$\left( \Delta _{\alpha ,\beta }\otimes id\right) \Phi
_{\alpha ,\beta \gamma }^{l}=\left( id\otimes \Phi _{\beta ,\gamma
}^{l}\right) \Phi _{\alpha \beta ,\gamma }^{l} $$
 Finally , for any $\alpha \in \pi , q\in A_{\alpha }^{2},q=a\otimes b$
\begin{eqnarray*}
\left( \varepsilon \otimes id\right) \Phi _{1,\alpha }^{l}\left(
q\right) &=&\left( \varepsilon \otimes id\right) \Phi _{1,\alpha }^{l}\left( a\otimes b\right) \\
&=&\left( \varepsilon \otimes id\right) \left( a_{\left(
1,1\right) }b_{\left( 1,1\right) }\otimes a_{\left( 2,\alpha
\right) }\otimes b_{\left( 2,\alpha \right) }\right) \\
&=&\varepsilon \left( a_{\left( 1,1\right) }b_{\left( 1,1\right)
}\right) a_{\left( 2,\alpha \right) }\otimes b_{\left( 2,\alpha
\right) }\\
&=&\varepsilon \left( a_{\left( 1,1\right) }\right) \varepsilon
\left( b_{\left( 1,1\right) }\right) a_{\left( 2,\alpha \right)
}\otimes b_{\left( 2,\alpha \right) }\\
 &=&a\otimes b\\
&=&q\\
\end{eqnarray*}
and thus the conditions of definition 3.1.1 are fulfilled and
$A^{2}=\left(
\left\{ A_{\alpha }^{2}\right\} _{\alpha \in \pi },\Phi ^{l}\right) $ is a $%
\pi -$ graded left covariant bimodule over $A .$ Similarly , one
can check that $A^{2}=\left( \left\{ A_{\alpha }^{2}\right\}
_{\alpha \in \pi },\Phi ^{r}\right) $ is a $\pi -$ graded right
covariant bimodule over $A.$  Finally , we check the bicovariance
condition

For any $\alpha ,\beta ,\gamma \in \pi , q\in A_{\alpha \beta
\gamma }^{2} , q=a\otimes b$ we compute
\begin{eqnarray*}
\left( id\otimes \Phi _{\beta ,\gamma }^{r}\right) \Phi _{\alpha
,\beta \gamma }^{l}\left( q\right)& =&\left( id\otimes \Phi
_{\beta ,\gamma }^{r}\right) \Phi
_{\alpha ,\beta \gamma }^{l}\left( a\otimes b\right) \\
&=&\left( id\otimes \Phi _{\beta ,\gamma }^{r}\right) \left(
a_{\left( 1,\alpha \right) }b_{\left( 1,\alpha \right) }\otimes
a_{\left( 2,\beta \gamma \right) }\otimes
b_{\left( 2,\beta \gamma \right) }\right) \\
&=&a_{\left( 1,\alpha \right) }b_{\left( 1,\alpha \right) }\otimes
a_{\left( 2,\beta \right) }\otimes b_{\left( 2,\beta \right)
}\otimes a_{\left( 3,\gamma \right) }b_{\left( 3,\gamma \right)
}\\
\end{eqnarray*}
\begin{eqnarray*}
\left( \Phi _{\alpha ,\beta }^{l}\otimes id\right) \Phi _{\alpha
\beta ,\gamma }^{r}\left( q\right) &=&\left( \Phi _{\alpha ,\beta
}^{l}\otimes id\right) \Phi _{\alpha
\beta ,\gamma }^{r}\left( a\otimes b\right) \\
&=&\left( \Phi _{\alpha ,\beta }^{l}\otimes id\right) \left(
a_{\left( 1,\alpha \beta \right) }\otimes b_{\left( 1,\alpha \beta
\right) }\otimes a_{\left( 2,\gamma
\right) }b_{\left( 2,\gamma \right) }\right)\\
&=&a_{\left( 1,\alpha \right) }b_{\left( 1,\alpha \right) }\otimes
a_{\left( 2,\beta \right) }\otimes b_{\left( 2,\beta \right)
}\otimes a_{\left( 3,\gamma \right) }b_{\left( 3,\gamma \right)
}\\
\end{eqnarray*}
which proves that $A^{2}=\left( \left\{ A_{\alpha }^{2}\right\}
_{\alpha \in \pi },\Phi ^{l},\Phi ^{r}\right) $ is a $\pi -$
graded bicovariant bimodule over $A.$

On $A\otimes A=\left\{ A_{\alpha }\otimes A_{\alpha }\right\}
_{\alpha \in \pi }$ we define two families of linear mappings

$r=\left\{ r_{\alpha }:A_{\alpha }\otimes A_{\alpha
}\longrightarrow A_{\alpha }\otimes A_{1}\right\} _{\alpha \in \pi
}$

 $t=\left\{ t_{\alpha }:A_{\alpha }\otimes A_{\alpha
}\longrightarrow A_{1}\otimes A_{\alpha }\right\} _{\alpha \in \pi
}$

For any $\alpha \in \pi ,a,b\in A_{\alpha }$ we set

$$r_{\alpha }\left( a\otimes b\right) =\left( a\otimes 1_{1}\right)
\Delta _{\alpha ,1}\left( b\right)\eqno(4.5)$$
$$t_{\alpha}\left( a\otimes b\right) =\left( 1_{1}\otimes a\right) \Delta
_{1,\alpha }\left( b\right)\eqno(4.6)$$  It is clear that
$r_{\alpha },t_{\alpha }$ are bijections for each $\alpha
\in \pi $ for example for $a\in A_{\alpha },b\in A_{1}$ the inverse of $%
r_{\alpha }$ is given by

$$r_{\alpha }^{-1}\left( a\otimes b\right) =\left( a\otimes
1_{\alpha }\right) \left( S_{\alpha ^{-1}}\otimes id\right) \Delta
_{\alpha ^{-1},\alpha }\left( b\right)\eqno(4.7)$$

Similarly , for $a\in A_{1},b\in A_{\alpha }$ the inverse of
$t_{\alpha }$ is given by

$$t_{\alpha }^{-1}\left( a\otimes b\right) =\left( b\otimes
1_{\alpha }\right) \left( S_{\alpha ^{-1}}\otimes id\right) \sigma
_{A_{\alpha ^{-1}},A_{\alpha }}\Delta _{\alpha ,\alpha
^{-1}}\left( b\right)\eqno(4.8)$$ One can easily show that for
each $\alpha \in \pi $ $,r_{\alpha }\left(
A_{\alpha }^{2}\right) =A_{\alpha }\otimes \ker \varepsilon $ , for let $%
\alpha \in \pi $ , $a\in A_{\alpha },b\in \ker \varepsilon $
\begin{eqnarray*}
m_{\alpha }r_{\alpha }^{-1}\left( a\otimes b\right) &=&m_{\alpha
}\left( \left( a\otimes 1_{\alpha }\right) \left( S_{\alpha
^{-1}}\otimes id\right) \Delta _{\alpha ^{-1},\alpha
}\left( b\right) \right) \\
&=&aS_{\alpha ^{-1}}\left( b_{\left( 1,\alpha ^{-1}\right)
}\right) b_{\left( 2,\alpha \right) }\\
&=&a\varepsilon \left( b\right) 1_{\alpha }\\
&=&0\\
\end{eqnarray*}
From which we get $r_{\alpha }^{-1}\left( A_{\alpha }\otimes \ker
\varepsilon \right) =A_{\alpha }^{2}$ i.e.
$$ r_{\alpha }\left(
A_{\alpha }^{2}\right) =A_{\alpha }\otimes \ker
\varepsilon\eqno(4.9)$$
 Similarly , one can prove that
$$t_{\alpha }\left( A_{\alpha }^{2}\right) =\ker \varepsilon
\otimes A_{\alpha }\eqno(4.10)$$

\begin{prop}
For any $\alpha ,\beta ,\gamma \in \pi$
$$ \left( \Delta _{\alpha
,\beta }\otimes id\right) r_{\alpha \beta }=\left( id\otimes
r_{\beta }\right) \Phi _{\alpha ,\beta }^{l}\eqno(4.11)$$

$$\left(id\otimes\Delta _{\alpha ,\beta }\right)t_{\alpha \beta
}=\left(t_{\alpha}\otimes id\right)\Phi _{\alpha ,\beta
}^{r}\eqno(4.12)$$

%==============================         =========================================
\end{prop}
\begin{proof}
We will prove that for any  $\alpha \in \pi $

$$r_{\alpha }=\left( id\otimes \varepsilon \otimes id \right) \Phi
_{\alpha ,1}^{l}\eqno(4.13)$$
$$t_{\alpha }=\left( \varepsilon \otimes id\otimes id\right) \Phi
_{1,\alpha }^{r}\eqno(4.14)$$

For any $\alpha \in \pi , a , b\in A_{\alpha } , a\otimes b\in
A_{\alpha }\otimes A_{\alpha }$

\begin{eqnarray*}
\left( id\otimes \varepsilon \otimes id\right) \Phi _{\alpha
,1}^{l}\left( a\otimes b\right) &=&\left( id\otimes \varepsilon
\otimes id\right) \left( a_{\left( 1,\alpha \right) }b_{\left(
1,\alpha \right) }\otimes a_{\left( 2,1\right) }\otimes b_{\left(
2,1\right) }\right) \\
&=&a_{\left( 1,\alpha \right) }b_{\left( 1,\alpha \right)
}\varepsilon \left( a_{\left( 2,1\right) }\right) \otimes
b_{\left( 2,1\right) }\\
&=&ab_{\left( 1,\alpha \right) }\otimes b_{\left( 2,1\right) }\\
&=&\left( a\otimes 1_{1}\right) \Delta _{\alpha ,1}\left(
b\right)\\
&=&r_{\alpha }\left( a\otimes b\right) \\
\end{eqnarray*}
\begin{eqnarray*}
\left( \varepsilon \otimes id\otimes id\right) \Phi _{1,\alpha
}^{r}\left( a\otimes b\right) &=&\left( \varepsilon \otimes
id\otimes id\right) \left( a_{\left( 1,1\right) }\otimes b_{\left(
1,1\right) }\otimes
a_{\left( 1,\alpha \right) }b_{\left( 1,\alpha \right) }\right) \\
&=&\varepsilon \left( a_{\left( 1,1\right) }\right) b_{\left(
1,1\right) }\otimes a_{\left( 1,\alpha \right) }b_{\left( 1,\alpha
\right) }\\
&=&b_{\left( 1,1\right) }\otimes ab_{\left( 1,\alpha \right) }\\
&=&\left( 1_{1}\otimes a\right) \Delta _{1,\alpha }\left(
b\right)\\
&=&t_{\alpha }\left( a\otimes b\right) \\
\end{eqnarray*}
To prove  4.11
\begin{eqnarray*}
\left( \Delta _{\alpha ,\beta }\otimes id\right) r_{\alpha \beta
}&=&\left( \Delta _{\alpha ,\beta }\otimes id\right) \left(
id\otimes \varepsilon \otimes
id\right) \Phi _{\alpha \beta ,1}^{l}\\
&=&\left( id\otimes id\otimes \varepsilon \otimes id\right) \left(
\Delta _{\alpha
,\beta }\otimes id\right) \Phi _{\alpha \beta ,1}^{l}\\
&=&\left( id\otimes id\otimes \varepsilon \otimes id\right) \left(
id\otimes \Phi _{\beta ,1}^{l}\right) \Phi _{\alpha ,\beta }^{l}\\
&=&\left( id\otimes \left( id\otimes \varepsilon \otimes id\right)
\Phi _{\beta ,1}^{l}\right)
\Phi _{\alpha ,\beta }^{l}\\
&=&\left( id\otimes r_{\beta }\right) \Phi _{\alpha
,\beta }^{l}\\
\end{eqnarray*}
Similarly ,one can prove 4.12.

\end{proof}

\begin{prop}

For any $\alpha \in \pi $ an element of $A_{\alpha }^{2}$ is left
-$\left( \textrm{right- respectively}\right) $invariant if and
only if it is of the form $r_{\alpha }^{-1}\left( 1_{\alpha }\otimes x\right) $ $%
\left( t_{\alpha }^{-1}\left( y\otimes 1_{\alpha }\right) \textrm{ respectively%
}\right) $where $x\in \ker \varepsilon $ $\left( y\in \ker
\varepsilon \textrm{ respectively}\right) .$
\end{prop}
\begin{proof}

For any $\alpha \in \pi $ , let $x\in \ker \varepsilon $ .We
compute
\begin{eqnarray*}
\Phi _{1,\alpha }^{l}\left( r_{\alpha }^{-1}\left( 1_{\alpha
}\otimes x\right) \right) &=&\Phi _{1,\alpha }^{l}\left( S_{\alpha
^{-1}}\left( x_{\left( 1,\alpha ^{-1}\right) }\right)\otimes \
x_{\left(
2,\alpha \right) }\right) \\
&=&S_{1}\left( x_{\left( 2,1\right) }\right) x_{\left( 3,1\right)
}\otimes S_{\alpha ^{-1}}\left( x_{\left( 1,\alpha ^{-1}\right)
}\right) \otimes x_{\left( 4,\alpha \right) }\\
&=&1_{1}\otimes \varepsilon \left( x_{\left( 2,1\right) }\right)
S_{\alpha ^{-1}}\left( x_{\left( 1,\alpha ^{-1}\right) }\right)
\otimes x_{\left( 3,\alpha \right) }\\
&=&1_{1}\otimes S_{\alpha ^{-1}}\left( x_{\left( 1,\alpha
^{-1}\right) }\right) \otimes x_{\left( 2,\alpha \right) }\\
&=&1_{1}\otimes r_{\alpha }^{-1}\left( 1_{\alpha }\otimes
x\right)\\
\end{eqnarray*}

i.e. $r_{\alpha }^{-1}\left( 1_{\alpha }\otimes x\right) $ is left
-invariant element.

Conversly,if $r_{\alpha }^{-1}\left( 1_{\alpha }\otimes x\right) $
is left -invariant element for some $\alpha \in \pi $ , let $x\in
\ker \varepsilon , $

$a\in A_{\alpha }.$ Equation 4.11 implies that $$\left( id\otimes
r_{\alpha }\right) \Phi _{1,\alpha }^{l}\left( r_{\alpha
}^{-1}\left( a \otimes x\right) \right) =\left( \Delta _{1,\alpha
}\otimes id\right) r_{\alpha }\left( r_{\alpha }^{-1}\left( a
\otimes x\right) \right) $$

From which we obtain
 $$1_{1}\otimes a\otimes x=\Delta _{1,\alpha }\left( a\right)
\otimes x$$ i.e.
 $$\Delta _{1,\alpha }\left( a\right) =1_{1}\otimes a$$

From which we obtain $$a=1_{\alpha }.$$
\end{proof}

\begin{thm}

Let $R$ be a right ideal of $A_{1}$ contained in $ker\varepsilon $ $%
,N=\left\{ N_{\alpha }\right\} _{\alpha \in \pi }$, where for each
$\alpha \in \pi ,$ $N_{\alpha }=$ $r_{\alpha }^{-1}\left(
A_{\alpha }\otimes R\right) $ is a sub-bimodule of $A^{2}=\left\{
A_{\alpha }^{2}\right\} _{\alpha \in \pi }.$ Moreover , let
$\Gamma =\left\{ \Gamma _{\alpha }\right\} _{\alpha \in \pi },$
$\Gamma _{\alpha }=A_{\alpha }^{2}/N_{\alpha },\Pi =\left\{ \Pi
_{\alpha }:A_{\alpha }^{2}\longrightarrow A_{\alpha
}^{2}/N_{\alpha }\right\} $ be the family of canonical epimorphisms, $%
d=\left\{ d_{\alpha }:d_{\alpha }=\Pi _{\alpha }\circ D_{\alpha }\right\} .$%
Then the $\pi -$graded first order differential calculus $\Gamma
=\left( \left\{ \Gamma _{\alpha }\right\} _{\alpha \in \pi
},d\right) $ is left covariant . Any $\pi -$graded left covariant
first order differential calculus on $A$ can be obtained in this
way.
\end{thm}
\begin{proof}

For any $\alpha \in \pi ,$let $R$ be a right ideal of $A_{1}$ contained in $%
ker\varepsilon .$ We shall prove that $r_{\alpha }^{-1}\left(
A_{\alpha }\otimes R\right) $ is a sub-bimodule of $A_{\alpha
}^{2}.$For any $\alpha
\in \pi ,$let $q\in r_{\alpha }^{-1}\left( A_{\alpha }\otimes R\right) ,$%
i.e. $q=r_{\alpha }^{-1}\left( b\otimes c\right) ,b\in A_{\alpha },c\in R.$%
For $a\in A_{\alpha }$
\begin{eqnarray*}
a\cdot q&=&\left( a\otimes 1_{\alpha }\right) q\\
&=&r_{\alpha}^{-1}\left( r_{\alpha }\left( \left( a\otimes
1_{\alpha }\right) q\right) \right) \\
&=&r_{\alpha }^{-1}\left( r_{\alpha }\left( a\otimes 1_{\alpha
}\right) r_{\alpha }\left( q\right) \right) \\
&=&r_{\alpha }^{-1}\left( \left( a\otimes 1_{1}\right) r_{\alpha
}\left( q\right) \right) \\
&=&r_{\alpha }^{-1}\left( \left( a\otimes 1_{1}\right) \left(
b\otimes c\right) \right) \\
&=&r_{\alpha }^{-1}\left( ab\otimes c\right) \\
&\in &r_{\alpha }^{-1}\left( A_{\alpha }\otimes R\right) \\
\end{eqnarray*}
\begin{eqnarray*}
 q\cdot a&=&q\left( 1_{\alpha }\otimes a\right) \\
&=&r_{\alpha }^{-1}\left( r_{\alpha }\left( q\left( 1_{\alpha
}\otimes a\right) \right) \right) \\
 &=&r_{\alpha }^{-1}\left( r_{\alpha }\left( q\right) \Delta
_{\alpha ,1}\left( a\right) \right) \\
 &=&r_{\alpha }^{-1}\left( \left( b\otimes c\right) \Delta
_{\alpha ,1}\left( a\right) \right) \\
&=&r_{\alpha }^{-1}\left( ba_{\left( 1,\alpha \right) }\otimes
ca_{\left( 2,1\right) }\right) \\
&\in &r_{\alpha }^{-1}\left( A_{\alpha }\otimes R\right) \\
\end{eqnarray*}
Which proves that $N_{\alpha }=r_{\alpha }^{-1}\left( A_{\alpha
}\otimes R\right) $ is a sub-bimodule of $A_{\alpha }^{2}.$

To prove that it is left covariant we have to prove that for any
$\alpha ,\beta \in \pi ,$ $\Phi _{\alpha ,\beta }^{l}\left(
N_{\alpha \beta }\right) \subset A_{\alpha }\otimes N_{\beta }.$

Using 4.11 we have
 $$\left( id\otimes r_{\beta
}\right) \Phi _{\alpha ,\beta }^{l}=\left( \Delta _{\alpha ,\beta
}\otimes id\right) r_{\alpha \beta }$$ i.e.
 $$\Phi _{\alpha,\beta }^{l}=\left( id\otimes r_{\beta }^{-1}\right)
\left( \Delta _{\alpha ,\beta }\otimes id\right) r_{\alpha \beta
}$$ Now ,for any $\alpha ,\beta \in \pi ,$ consider $N_{\alpha
\beta }=r_{\alpha \beta }^{-1}\left( A_{\alpha \beta }\otimes
R\right) $
\begin{eqnarray*}
\Phi _{\alpha ,\beta }^{l}\left( N_{\alpha \beta }\right)
&=&\left( id\otimes r_{\beta }^{-1}\right) \left( \Delta _{\alpha
,\beta }\otimes id\right) r_{\alpha \beta
}\left( N_{\alpha \beta }\right) \\
&=&\left( id\otimes r_{\beta }^{-1}\right) \left( \Delta _{\alpha
,\beta }\otimes id\right) r_{\alpha \beta } \left( r_{\alpha \beta
}^{-1}\left( A_{\alpha \beta }\otimes
R\right) \right) \\
&=&\left( id\otimes r_{\beta }^{-1}\right) \left( \Delta _{\alpha
,\beta }\otimes id\right) \left( A_{\alpha
\beta }\otimes R\right) \\
\end{eqnarray*}

\begin{eqnarray*}
&=&\left( id\otimes r_{\beta }^{-1}\right) \left( \Delta _{\alpha
,\beta }\left( A_{\alpha \beta }\right) \otimes
R\right) \\
&\subset& \left( id\otimes r_{\beta }^{-1}\right)
\left( A_{\alpha }\otimes A_{\beta }\otimes R\right) \\
&=&A_{\alpha }\otimes r_{\beta }^{-1}\left( A_{\beta }\otimes
R\right) \\
&=&A_{\alpha }\otimes N_{\beta }\\
\end{eqnarray*}
Conversly , if $N=\left( \left\{ N_{\alpha }\right\} _{\alpha \in
\pi },\Phi ^{l}\right) $ is a left covariant bimodule , then ,
using theorem 3.3.1 and proposition 4.2 there exsists a family
$\left( x_{i}\right) _{i\in I}$ of elements of $ker\varepsilon $
such that for any $\alpha \in \pi ,q\in N_{\alpha }$ can be
written as $q=\sum_{i}a_{i}\cdot $ $r_{\alpha }^{-1}\left(
1_{\alpha }\otimes x_{i}\right) ,a_{i}\in A_{\alpha }.$ But for
each $i\in I$ we have
\begin{eqnarray*}
a_{i}\cdot r_{\alpha }^{-1}\left( 1_{\alpha }\otimes
x_{i}\right)&=&\left( a_{i}\otimes 1_{\alpha }\right) r_{\alpha
}^{-1}\left( 1_{\alpha }\otimes x_{i}\right) \\
&=&r_{\alpha }^{-1}\left( r_{\alpha }\left( a_{i}\otimes 1_{\alpha
}\right) \left( 1_{\alpha }\otimes x_{i}\right) \right) \\
&=&r_{\alpha }^{-1}\left( \left( a_{i}\otimes 1_{\alpha }\right)
\left( 1_{\alpha }\otimes x_{i}\right) \right) \\
&=&r_{\alpha }^{-1}\left( a_{i}\otimes x_{i}\right) \\
\end{eqnarray*}
 Denoting by $R_{\alpha }$ the linear span of all $x_{i}^{,}s$
we obtain that
 $N_{\alpha }=r_{\alpha }^{-1}\left( A_{\alpha }\otimes R_{\alpha }\right) $

We shall show that all $R_{\alpha }^{,}s$ coincide with $R_{1}.$
From proposition 4.2 we have

 $$_{inv}N_{\alpha }=r_{\alpha }^{-1}\left( 1_{\alpha }\otimes R_{\alpha }\right) $$
and since $N_{\alpha }$ is a left covariant bimodule we have
\begin{eqnarray*}
\Phi _{\alpha ,1}^{l}\left( _{inv}N_{\alpha }\right) &=&1_{\alpha
}\otimes _{inv}N_{1}\\
 &=&1_{\alpha }\otimes
r_{1}^{-1}\left(1_{1}\otimes R_{1}\right)\\
\end{eqnarray*}
Now let $r_{\alpha }^{-1}\left( 1_{\alpha }\otimes x_{i}\right)
\in _{inv}N_{\alpha }$ , $x_{i}\in R_{\alpha }$
\begin{eqnarray*}
\Phi _{\alpha ,1}^{l}\left( r_{\alpha }^{-1}\left( 1_{\alpha
}\otimes x_{i}\right) \right)&=&\Phi _{\alpha ,1}^{l}\left(
S_{\alpha ^{-1}}\left( x_{i\left( 1,\alpha ^{-1}\right) }\right)
\otimes x_{i\left( 2,\alpha \right) }\right) \\
&=&S_{\alpha ^{-1}}\left( x_{i\left( 2,\alpha ^{-1}\right)
}\right) x_{i\left( 3,\alpha \right) }\otimes S_{1^{-1}}\left(
x_{i\left( 1,1\right) }\right) \otimes x_{i\left( 4,1\right) }\\
&=&1_{\alpha }\varepsilon \left( x_{i\left( 2,1\right) }\right)
\otimes S_{1^{-1}}\left( x_{i\left( 1,1\right) }\right) \otimes
x_{i\left( 3,1\right) }\\
&=&1_{\alpha }\otimes S_{1^{-1}}\left( x_{i\left( 1,1\right)
}\right) \otimes x_{i\left( 2,1\right) }\\
&=&1_{\alpha }\otimes r_{1}^{-1}\left( 1_{1}\otimes x_{i}\right)\\
\end{eqnarray*}
i.e. $$x_{i}\in R_{1}$$ $$\Longrightarrow R_{\alpha }\subseteq
R_{1}$$.\\ Similarly we can show that $R_{1}\subseteq R_{\alpha
},$ and hence $R_{\alpha }=R_{1}$ for each $\alpha \in \pi .$
Denote by $R$ to any of the $R_{\alpha }^{,}s ,$ then  $$N_{\alpha
}= r_{\alpha}^{-1}\left( A_{\alpha }\otimes R\right) $$\\ It
remains to show that $R$ is a right ideal of $A_{1}.$ Let $x\in R
, a\in A_{1} , $ then $r_{1}^{-1}\left( 1_{1}\otimes x\right) \in
N_{1}.$
\begin{eqnarray*}
 r_{1}^{-1}\left( 1_{1}\otimes
x\right) \cdot a&=&r_{1}^{-1}\left( 1_{1}\otimes x\right) \left(
1_{1}\otimes a\right) \\
&=&r_{1}^{-1}\left( \left( 1_{1}\otimes x\right) r_{1}\left(
1_{1}\otimes a\right) \right)\\
&\in& N_{1}=r_{1}^{-1}\left(A_{1}\otimes R\right) \\
&&\left( N_{1}\textrm{ is a bimodule}\right)\\
\end{eqnarray*}

i.e $$\left( 1_{1}\otimes x\right) r_{1}\left( 1_{1}\otimes
a\right) \in A_{1}\otimes R$$ therefore
\begin{eqnarray*}
 \left( 1_{1}\otimes x\right) r_{1}\left( 1_{1}\otimes
a\right) &=&r_{1}\left( r_{1}^{-1}\left( 1_{1}\otimes x\right)
\left( 1_{1}\otimes a\right) \right) \\
&=&r_{1}\left( r_{1}^{-1}\left( \left( 1_{1}\otimes x\right)
\Delta _{1,1}\left( a\right) \right) \right) \\
&=&\left( 1_{1}\otimes x\right) \Delta _{1,1}\left( a\right) \in
A_{1}\otimes R\\
\end{eqnarray*}
and $\left( \varepsilon \otimes id\right) \left( \left(
1_{1}\otimes x\right) \Delta _{1,1}\left( a\right) \right) =xa\in
R$
\end{proof}

\begin{thm}

Let $R$ be a right ideal of $A_{1}$ contained in $ker\varepsilon $ $%
,N=\left\{ N_{\alpha }\right\} _{\alpha \in \pi }$,where for each
$\alpha \in \pi ,$ $N_{\alpha }=$ $t_{\alpha }^{-1}\left(
A_{\alpha }\otimes R\right) $ is a sub-bimodule of $A^{2}=\left\{
A_{\alpha }^{2}\right\} _{\alpha \in \pi }.$ Moreover ,let $\Gamma
=\left\{ \Gamma _{\alpha }\right\} _{\alpha \in \pi },$ $\Gamma
_{\alpha }=A_{\alpha }^{2}/N_{\alpha },\Pi =\left\{ \Pi _{\alpha
}:A_{\alpha }^{2}\longrightarrow A_{\alpha
}^{2}/N_{\alpha }\right\} $ be the family of canonical epimorphisms , $%
d=\left\{ d_{\alpha }:d_{\alpha }=\Pi _{\alpha }\circ D_{\alpha }\right\} .$%
Then the first order differential calculus $\Gamma =\left( \left\{
\Gamma _{\alpha }\right\} _{\alpha \in \pi },d\right) $ is right
covariant . Any right covariant first order differential calculus
on $A$ can be obtained in this way.
\end{thm}
\begin{proof}
For any $\alpha \in \pi ,$let $R$ be a right ideal of $A_{1}$ contained in $%
ker\varepsilon .$ We shall prove that $t_{\alpha }^{-1}\left(
R\otimes A_{\alpha }\right) $ is a sub-bimodule of $A_{\alpha
}^{2}.$For any $\alpha
\in \pi ,$let $q\in t_{\alpha }^{-1}\left( R\otimes A_{\alpha }\right) ,$%
i.e. $q=t_{\alpha }^{-1}\left( d\otimes e\right) ,d\in R,e\in A_{\alpha }.$%
For $a\in A_{\alpha }$
\begin{eqnarray*}
a\cdot q&=&\left( a\otimes 1_{\alpha }\right) q\\
&=&t_{\alpha }^{-1}\left( t_{\alpha }\left( \left( a\otimes
1_{\alpha }\right) q\right) \right) \\
&=&t_{\alpha }^{-1}\left( t_{\alpha }\left( a\otimes 1_{\alpha
}\right) t_{\alpha }\left( q\right) \right) \\
&=&t_{\alpha }^{-1}\left( \left( 1_{1}\otimes a\right) t_{\alpha
}\left( q\right) \right) \\
&=&t_{\alpha }^{-1}\left( \left( 1_{1}\otimes a\right) \left(
d\otimes e\right) \right) \\
&=&t_{\alpha }^{-1}\left( d\otimes ae\right) \\
&\in& t_{\alpha }^{-1}\left( R\otimes A_{\alpha }\right) \\
\end{eqnarray*}
\begin{eqnarray*}
q\cdot a&=&q\left( 1_{\alpha }\otimes a\right) \\
&=&t_{\alpha }^{-1}\left( t_{\alpha }\left( q\left( 1_{\alpha
}\otimes a\right) \right) \right) \\
&=&t_{\alpha }^{-1}\left( t_{\alpha }\left( q\right) \Delta
_{1,\alpha }\left( a\right) \right) \\
&=&r_{\alpha }^{-1}\left( \left( d\otimes e\right) \Delta
_{1,\alpha }\left( a\right) \right) \\
&=&t_{\alpha }^{-1}\left( da_{\left( 1,1\right) }\otimes
ca_{\left( 2,\alpha \right) }\right) \\
&\in& t_{\alpha }^{-1}\left( R\otimes A_{\alpha }\right) \\
\end{eqnarray*}
which proves that $N_{\alpha }=t_{\alpha }^{-1}\left( R\otimes
A_{\alpha }\right) $ is a sub-bimodule of $A_{\alpha }^{2}.$
 To prove that it is right covariant we have to prove that for any
$\alpha ,\beta \in \pi ,$ $\Phi _{\alpha ,\beta }^{l}\left(
N_{\alpha \beta }\right) \subset N_{\alpha }\otimes A_{\beta }.$
Using 4.13 we have
 $$\left( t_{\alpha }\otimes id\right) \Phi _{\alpha ,\beta }^{r}=\left( id\otimes
\Delta _{\alpha ,\beta }\right) t_{\alpha \beta }$$
 i.e.
$$\Phi _{\alpha ,\beta }^{r}=\left( t_{\alpha }^{-1}\otimes
id\right) \left( id\otimes \Delta _{\alpha ,\beta }\right)
t_{\alpha \beta }$$ Now ,for any $\alpha ,\beta \in \pi ,$consider
$N_{\alpha \beta }=t_{\alpha \beta }^{-1}\left( R\otimes A_{\alpha
\beta }\right) $
\begin{eqnarray*}
\Phi _{\alpha ,\beta }^{l}\left( N_{\alpha \beta }\right)
&=&\left( t_{\alpha }^{-1}\otimes id\right) \left( id\otimes
\Delta _{\alpha ,\beta }\right) t_{\alpha \beta
}\left( N_{\alpha \beta }\right) \\
&=&\left( t_{\alpha }^{-1}\otimes id\right) \left( id\otimes
\Delta _{\alpha ,\beta }\right) t_{\alpha \beta }\left( t_{\alpha
\beta }^{-1}\left( R\otimes A_{\alpha \beta
}\right) \right) \\
&=&\left( t_{\alpha }^{-1}\otimes id\right) \left( id\otimes
\Delta _{\alpha ,\beta }\right) \left( R\otimes
A_{\alpha \beta }\right) \\
&\subset &\left( t_{\alpha }^{-1}\otimes id\right)
\left( R\otimes A_{\alpha }\otimes A_{\beta }\right) \\
&=&t_{\alpha }^{-1}\left( R\otimes A_{\alpha }\right) \otimes
A_{\beta }\\
&=&N_{\alpha }\otimes A_{\beta }\\
\end{eqnarray*}
Conversly , if $N=\left( \left\{ N_{\alpha }\right\} _{\alpha \in
\pi },\Phi ^{r}\right) $ is a right covariant bimodule then ,
using theorem 3.4.1 and proposition 4.2 there exsists a family
$\left( y_{i}\right) _{i\in I}$ of elements of $ker\varepsilon $
such that for any $\alpha \in \pi ,q\in N_{\alpha }$ can be
written as $q=\sum_{i}a_{i}\cdot $ $t_{\alpha }^{-1}\left(
y_{i}\otimes 1_{\alpha }\right) ,a_{i}\in A_{\alpha }.$ But for
each $i\in I$ we have
\begin{eqnarray*}
a_{i}\cdot r_{\alpha }^{-1}\left( y_{i}\otimes 1_{\alpha }\right)
&=&\left( a_{i}\otimes 1_{\alpha }\right) t_{\alpha
}^{-1}\left( y_{i}\otimes 1_{\alpha }\right) \\
&=&t_{\alpha }^{-1}\left( t_{\alpha }\left( 1_{\alpha }\otimes
a_{i}\right) \left( y_{i}\otimes 1_{\alpha }\right) \right)\\
&=&t_{\alpha }^{-1}\left( \left( 1_{\alpha }\otimes a_{i}\right)
\left( y_{i}\otimes 1_{\alpha }\right) \right) \\
&=&t_{\alpha }^{-1}\left( y_{i}\otimes a_{i}\right) \\
\end{eqnarray*}
Denoting by $R_{\alpha }$ the linear span of all $x_{i}^{,}s$ we
obtain that
 $$N_{\alpha }=t_{\alpha}^{-1}\left( R\otimes A_{\alpha }\right)
 $$
We shall show that all $R_{\alpha }^{,}s$ coincide with $R_{1} .$
From proposition 4.2 we have $$N_{\alpha }^{inv}=t_{\alpha
}^{-1}\left( 1_{\alpha }\otimes A_{\alpha }\right) $$

and since $N_{\alpha }$ is a left covariant bimodule we have
\begin{eqnarray*}
\Phi _{1,\alpha }^{r}\left( N_{\alpha }^{inv}\right)
&=&N_{1}^{inv}\otimes 1_{\alpha }\\
&=&t_{1}^{-1}\left( R\otimes A_{1}\right) \otimes 1_{\alpha}\\
\end{eqnarray*}
 Now let $t_{\alpha }^{-1}\left( y_{i}\otimes
1_{\alpha }\right) \in N_{\alpha }^{inv}$ , $y_{i}\in R_{\alpha }$
\begin{eqnarray*}
\Phi _{1,\alpha }^{r}\left( t_{\alpha }^{-1}\left( y_{i}\otimes
1_{\alpha }\right) \right) &=&\Phi _{1,\alpha }^{r}\left(
S_{\alpha ^{-1}}\left( y_{i\left( 2,\alpha ^{-1}\right) }\right)
\otimes y_{i\left( 1,1\right) }\right) \\
&=&S_{1^{-1}}\left( y_{i\left( 4,1\right) }\right) \otimes
y_{i\left( 1,1\right) }\otimes S_{\alpha ^{-1}}\left( y_{i\left(
3,\alpha ^{-1}\right) }\right) y_{i\left( 2,\alpha \right) }\\
&=&S_{1^{-1}}\left( y_{i\left( 3,1\right) }\right) \otimes
y_{i\left( 1,1\right) }\otimes 1_{\alpha }\varepsilon \left(
y_{i\left( 2,1\right) }\right) \\
&=&S_{1^{-1}}\left( y_{i\left( 2,1\right) }\right) \otimes
y_{i\left( 1,1\right) }\otimes 1_{\alpha }\\
&=&t_{1}^{-1}\left( y_{i}\otimes 1_{1}\right) \otimes 1_{\alpha}\\
\end{eqnarray*}

i.e. $$y_{i}\in R_{1}$$ $$\Longrightarrow R_{\alpha }\subseteq
R_{1}$$

Similarly we can show that $R_{1}\subseteq R_{\alpha },$ and hence
$ R_{\alpha }=R_{1}$ for each $\alpha \in \pi .$

Denote by $R$ to any of the $R_{\alpha }^{,}s,$ then
$N_{\alpha}=t_{\alpha }^{-1}\left( R\otimes A_{\alpha }\right) $

It remains to show that $R$ is a right ideal of $A_{1}.$

Let $y\in R,a\in A_{1},$ $A$ then $t_{1}^{-1}\left( y\otimes
1_{1}\right) \in N_{1}.$
\begin{eqnarray*}
t_{1}^{-1}\left( y\otimes 1_{1}\right) \cdot a&=&t_{1}^{-1}\left(
y\otimes 1_{1}\right) \left( 1_{1}\otimes a\right) \\
&=&t_{1}^{-1}\left( \left( y\otimes 1_{1}\right) t_{1}\left(
1_{1}\otimes a\right) \right) \in N_{1}\\
&=&t_{\alpha }^{-1}\left(R\otimes A_{1}\right) \\
&&\left( N_{1}\textrm{ is a bimodule}\right) \\
\end{eqnarray*}
i.e   $\left( y\otimes 1_{1}\right) t_{1}\left( 1_{1}\otimes
a\right) \in $ $R\otimes A_{1}$

therefore
\begin {eqnarray*}
\left( y\otimes 1_{1}\right) t_{1}\left( 1_{1}\otimes a\right)
&=&t_{1}\left( t_{1}^{-1}\left( y\otimes 1_{1}\right) \left(
1_{1}\otimes a\right) \right) \\
&=&t_{1}\left( t_{1}^{-1}\left( \left( y\otimes 1_{1}\right)
\Delta _{1,1}\left( a\right) \right) \right) \\
&=&\left( y\otimes 1_{1}\right) \Delta _{1,1}\left( a\right) \in
A_{1}\otimes R\\
\end{eqnarray*}
and $\left( id\otimes \varepsilon \right) \left( \left( y\otimes
1_{1}\right) \Delta _{1,1}\left( a\right) \right) =ya\in R$

\end{proof}
 We shall now formulate the concept of $ad-$invariance . Let
$$ad_{\alpha }:A_{1}\longrightarrow A_{1}\otimes A_{\alpha }$$ be
such that for any $a\in A_{1}$

$$ad_{\alpha }\left( a\right) =t_{\alpha }\left( r_{\alpha
}^{-1}\left( 1_{\alpha }\otimes a\right) \right)\eqno(4.15)$$ i.e.
$$ad_{\alpha }\left( a\right) =a_{\left( 2,1\right) }\otimes
S_{\alpha ^{-1}}\left( a_{\left( 1,\alpha ^{-1}\right)
}\right)a_{\left( 3,\alpha \right) }$$

where

$$\left( id\otimes \Delta _{1,\alpha }\right) \Delta _{\alpha
^{-1},\alpha }\left( a\right) = a_{\left( 1,\alpha ^{-1}\right)
}\otimes a_{\left( 2,1\right) }\otimes a_{\left( 3,\alpha \right)
}\eqno(4.16)$$
 such that

$$\left( ad_{\alpha }\otimes id\right) ad_{\beta }\left( a\right) =
\left( id\otimes \Delta _{\alpha ,\beta }\right) ad_{\alpha \beta
}\eqno(4.17)$$

 Using 4.15 ,and the standared properties of
comultiplication and coinverse one can prove 4.17, for let $a\in
A_{1}. $For any $\alpha ,\beta \in \pi $
\begin{eqnarray*}
\left( ad_{\alpha }\otimes id\right) ad_{\beta }\left( a\right)
&=&\left( ad_{\alpha }\otimes id\right) \left( a_{\left(
2,1\right) }\otimes S_{\beta ^{-1}}\left( a_{\left( 1,\beta
^{-1}\right) }\right) a_{\left( 3,\beta
\right) }\right) \\
&=&a_{\left( 3,1\right) }\otimes S_{\alpha ^{-1}}\left( a_{\left(
2,\alpha ^{-1}\right) }\right) a_{\left( 4,\alpha \right) }\otimes
S_{\beta ^{-1}}\left( a_{\left( 1,\beta ^{-1}\right) }\right)
a_{\left( 5,\beta \right) }\\
\end{eqnarray*}
\begin{eqnarray*}
\left( id\otimes \Delta _{\alpha ,\beta }\right) ad_{\alpha \beta
}&=&\left( id\otimes \Delta _{\alpha ,\beta }\right) \left(
a_{\left( 2,1\right) }\otimes S_{\left( \alpha \beta \right)
^{-1}}\left( a_{\left( 1,\left( \alpha \beta \right) ^{-1}\right)
}\right) a_{\left( 3,\alpha \beta
\right) }\right) \\
&=&a_{\left( 3,1\right) }\otimes S_{\alpha ^{-1}}\left( a_{\left(
2,\alpha ^{-1}\right) }\right) a_{\left( 4,\alpha \right) }\otimes
S_{\beta ^{-1}}\left( a_{\left( 1,\beta ^{-1}\right) }\right)
a_{\left( 5,\beta \right) }\\
\end{eqnarray*}
which proves equation 4.17.
\end{proof}
A linear subset $T\subset A_{1}$ is $\pi -ad$ invariant if
$ad_{\alpha }\left( T\right) \subset T\otimes A$ for any $\alpha
\in \pi .$

\begin{lem}

Let $T$ be $\pi -ad$ invariant subset of $A_{1},R$ be a right ideal of $%
A_{1} $ generated by $T$. Then $R$ is $\pi -ad$ invariant.
\end{lem}

\begin{proof}

Let $a,b\in A_{1},$we will prove that for any $\alpha \in \pi $

$$ad_{\alpha }\left( ab\right) =\left( 1_{1}\otimes S_{\alpha
^{-1}}\left( b_{\left( 1,\alpha ^{-1}\right) }\right) \right)
ad_{\alpha }\left( a\right) \Delta _{1,\alpha }\left( b_{\left(
2,\alpha \right) }\right)\eqno(4.18)$$
\begin{eqnarray*}
r_{\alpha }^{-1}\left( 1_{\alpha }\otimes ab\right)&=&\left(
1_{\alpha }\otimes 1_{\alpha }\right) \left( S_{\alpha
^{-1}}\otimes id\right) \Delta _{\alpha ^{-1},\alpha
}\left( ab\right) \\
&=&S_{\alpha ^{-1}}\left( a_{\left( 1,\alpha ^{-1}\right)
}b_{\left( 1,\alpha ^{-1}\right) }\right) \otimes a_{\left(
2,\alpha \right) }b_{\left( 2,\alpha \right) }\\
&=&S_{\alpha ^{-1}}\left( b_{\left( 1,\alpha ^{-1}\right) }\right)
S_{\alpha ^{-1}}\left( a_{\left( 1,\alpha ^{-1}\right) }\right)
\otimes a_{\left( 2,\alpha \right) }b_{\left( 2,\alpha \right) }\\
&=&\left( S_{\alpha ^{-1}}\left( b_{\left( 1,\alpha ^{-1}\right)
}\right) \otimes 1_{\alpha }\right) \left( S_{\alpha ^{-1}}\left(
a_{\left( 1,\alpha ^{-1}\right) }\right) \otimes a_{\left(
2,\alpha \right) }\right) \left( 1_{\alpha }\otimes b_{\left(
2,\alpha \right) }\right) \\
&=&\left( S_{\alpha ^{-1}}\left( b_{\left( 1,\alpha ^{-1}\right)
}\right) \otimes 1_{\alpha }\right) r_{\alpha }^{-1}\left(
a\right) \left( 1_{\alpha }\otimes b_{\left( 2,\alpha \right)
}\right) \\
\end{eqnarray*}
Applying $t_{\alpha }$ to both sides of the above equation we get
$$t_{\alpha }r_{\alpha }^{-1}\left( 1_{\alpha }\otimes ab\right)
=t_{\alpha }\left( S_{\alpha ^{-1}}\left( b_{\left( 1,\alpha
^{-1}\right) }\right) \otimes 1_{\alpha }\right) t_{\alpha
}r_{\alpha }^{-1}\left( a\right) t_{\alpha }\left( 1_{\alpha
}\otimes b_{\left( 2,\alpha \right) }\right)\eqno(4.19)$$
$$ad_{\alpha }\left( ab\right) =\left( 1_{\alpha }\otimes S_{\alpha
^{-1}}\left( b_{\left( 1,\alpha ^{-1}\right) }\right) \right)
ad_{\alpha }\left( a\right) \Delta _{1,\alpha }\left( b_{\left(
2,\alpha \right) }\right)\eqno(4.20)$$

Thus for $a,b\in T,a,b,ab\in R , R$ being an ideal in $A_{1},T$
being $\pi -ad$ invariant we find that

$$ad_{\alpha }\left( ab\right) \in R\otimes A_{\alpha }$$ i.e.
$$ad_{\alpha }\left( R\right) \subset R\otimes A_{\alpha }$$  which
means that $R$ is $\pi -ad$ invariant.
\end{proof}
%=====================================================================
Let $A^{2}=\left( \left\{ A_{\alpha }^{2}\right\} _{\alpha \in \pi
},\Phi ^{l},\Phi ^{r}\right) $ is a $\pi -$ graded bicovariant
bimodule over $A.$By virtue of condition 3 of definition 3.3 we
have for any $\alpha ,\beta ,\gamma \in \pi $
 $\left( \Phi_{\alpha ,1}^{l}\otimes id\right) \Phi _{\alpha
,\beta }^{r}=\left( id\otimes \Phi _{1,\beta }^{r}\right) \Phi
_{\alpha ,\beta }^{l}$

Applying $id\otimes \varepsilon \otimes id\otimes id $ to both
sides of the above equation we get

$$\left( \left( id\otimes \varepsilon \otimes id\right) \Phi
_{\alpha ,1}^{l}\otimes id\right) \Phi _{\alpha ,\beta
}^{r}=\left( id\otimes \left( \varepsilon \otimes id\otimes
id\right) \Phi _{1,\beta }^{r}\right) \Phi _{\alpha ,\beta }^{l}$$

Using 4.13, 4.14
 $$\left( r_{\alpha }\otimes id\right)
\Phi _{\alpha ,\beta }^{r}=\left( id\otimes t_{\beta }\right) \Phi
_{\alpha ,\beta }^{l}$$

i.e.  $$\Phi _{\alpha ,\beta }^{r}=\left( r_{\alpha }^{-1}\otimes
id\right) \left( id\otimes t_{\beta }\right) \Phi _{\alpha ,\beta
}^{l}$$

Now let $x\in \ker \varepsilon .$ From proposition 4.2 for any
$\alpha \in \pi $ we have $r_{\alpha }^{-1}\left( 1_{\alpha
}\otimes x\right) $ is a left invariant element then
\begin{eqnarray*}
\Phi _{\alpha ,\beta }^{r}\left( r_{\alpha \beta }^{-1}\left(
1_{\alpha \beta }\otimes x\right) \right) &=&\left( r_{\alpha
}^{-1}\otimes id\right) \left( id\otimes t_{\beta }\right) \Phi
_{\alpha ,\beta }^{l}\left( r_{\alpha \beta }^{-1}\left( 1_{\alpha
\beta }\otimes x\right)
\right) \\
&=&\left( r_{\alpha }^{-1}\otimes id\right) \left( id\otimes
t_{\beta }\right) \left( 1_{\alpha }\otimes r_{\beta }^{-1}\left(
1_{\beta }\otimes x\right)
\right)\\
&=&\left( r_{\alpha }^{-1}\otimes id\right) \left( 1_{\alpha
}\otimes t_{\beta }r_{\beta }^{-1}\left( 1_{\beta
}\otimes x\right) \right) \\
&=&\left( r_{\alpha }^{-1}\otimes id\right) \left(
1_{\alpha }\otimes ad_{\beta }\left( x\right) \right) \\
\end{eqnarray*}

\begin{thm}

Let $R$ be a right ideal of $A_{1}$ contained in $\ker \varepsilon $ and $%
\Gamma =\left( \left\{ \Gamma _{\alpha }\right\} _{\alpha \in \pi
},d\right) $be the $\pi -$graded left covariant first order
differential calculus described in 4.3 .Then $\Gamma =\left(
\left\{ \Gamma _{\alpha }\right\} _{\alpha \in \pi },d\right) $ is
bicovariant if and only if $R$ is $\pi -$ad invariant.
\end{thm}
\begin{proof}

Let for any $\alpha \in \pi $ $R$ be a right ideal of $A_{1}$such that $%
R\subset \ker \varepsilon $ and $N_{\alpha }=r_{\alpha
}^{-1}\left( A_{\alpha }\otimes R\right) .$Using theorem 4.3 we
see that $N=\left( \left\{ N_{\alpha }\right\} _{\alpha \in \pi
},\Phi ^{l}\right) $ is a $\pi - $graded left covariant bimodule
.Assume that $R$ is $\pi -$ad invariant ,
let $_{inv}N_{\alpha }$ be the set of all left invariant elements of $%
N_{\alpha }$ for each $\alpha \in \pi .$ Then formula 3.78 shows
that for any $\alpha ,\beta \in \pi $

$$\Phi _{\alpha ,\beta }^{r}\left( _{inv}N_{\alpha \beta }\right)
\subset _{inv}N_{\alpha }\otimes A_{\beta }$$

Now decomposition 3.26 shows that $\Phi _{\alpha ,\beta
}^{r}\left( N_{\alpha \beta }\right) \subset N_{\alpha }\otimes
A_{\beta },$and this means that implication 2.14 holds.

Conversly,assume that $N=\left\{ N_{\alpha }\right\} _{\alpha \in
\pi }$ is a $\pi -$graded bicovariant bimodule .This means that
2.14 holds.Then (see proof 4.4) for each $\alpha \in \pi ,$
$N_{\alpha }=t_{\alpha }^{-1}\left( R^{\prime }\otimes A_{\alpha
}\right) $ where $R^{\prime }$ be a right ideal
of $A_{1}$ such that $R^{\prime }\subset \ker \varepsilon .$In particular, $%
N_{1}=t_{1}^{-1}\left( R^{\prime }\otimes A_{1}\right) $ .Using
3.78 and that $\left( \varepsilon \otimes id\right)
t_{1}^{-1}\left( a\otimes b\right) =a\varepsilon \left( b\right) $
,and $\left( id\otimes \varepsilon \right) r_{1}^{-1}\left(
a\otimes b\right) =aS_{1}\left( b\right) ,$one can easily checks
that $R=R^{\prime }.$

So we have for any $\alpha \in \pi $

$$r_{\alpha }^{-1}\left( A_{\alpha }\otimes R\right) =t_{\alpha
}^{-1}\left( R\otimes A_{\alpha }\right) $$

$$t_{\alpha }r_{\alpha }^{-1}\left( A_{\alpha }\otimes R\right)
=R\otimes A_{\alpha }$$

therefore  $ad_{\alpha }\left( R\right) =t_{\alpha }r_{\alpha
}^{-1}\left( 1_{\alpha }\otimes R\right) $

$\subset t_{\alpha }r_{\alpha }^{-1}\left( A_{\alpha }\otimes
R\right) $
 $=R\otimes A_{\alpha }$

 therefore $R$ is $\pi-$ad invariant.
\end{proof}

\end{document}